\magnification=1200

\hsize=14cm    
\vsize=20.5cm       

\parindent=0cm   \parskip=0pt     
\pageno=1 

\def\ind{\hskip 1cm\relax}


\hoffset=15mm    
\voffset=1cm    
 

\ifnum\mag=\magstep1
\hoffset=-0.5cm   
\voffset=-0.5cm   
\fi

\pretolerance=500 \tolerance=1000  \brokenpenalty=5000

\catcode`\@=11


\newcount\secno
\newcount\prmno
\newif\ifnotfound
\newif\iffound

\def\namedef#1{\expandafter\def\csname #1\endcsname}
\def\nameuse#1{\csname #1\endcsname}

\long\def\ifundefined#1#2#3{\expandafter\ifx\csname
  #1\endcsname\relax#2\else#3\fi}
\def\hwrite#1#2{{\let\the=0\edef\next{\write#1{#2}}\next}}

\toksdef\ta=0 \toksdef\tb=2
\long\def\leftappenditem#1\to#2{\ta={\\{#1}}\tb=\expandafter{#2}%
                                \edef#2{\the\ta\the\tb}}
\long\def\rightappenditem#1\to#2{\ta={\\{#1}}\tb=\expandafter{#2}%
                                \edef#2{\the\tb\the\ta}}

\def\lop#1\to#2{\expandafter\lopoff#1\lopoff#1#2}
\long\def\lopoff\\#1#2\lopoff#3#4{\def#4{#1}\def#3{#2}}

\def\ismember#1\of#2{\foundfalse{\let\given=#1%
    \def\\##1{\def\next{##1}%
    \ifx\next\given{\global\foundtrue}\fi}#2}}

\def\section#1{\medbreak
               \global\def\currenvir{section}
               \global\advance\secno by1\global\prmno=0
               {\bf \number\secno. {#1}}}

\def\subsection{\global\def\currenvir{subsection}
                \global\advance\prmno by1
                {\hskip-.5truemm(\number\secno.\number\prmno)}\hskip 5truemm}

\def\formule{\global\def\currenvir{formule}
                \global\advance\prmno by1
                {\hbox{\rm (\number\secno.\number\prmno)}}}

\def\proclaim#1{\global\advance\prmno by 1
                {\bf #1 \the\secno.\the\prmno.-- }}


\def\ex#1{\medbreak\global\advance\prmno by 1
                {\bf #1 \the\secno.\the\prmno.}}

\def\qq#1{\medbreak\global\advance\prmno by 1
                { #1 \the\prmno)}}

\long\def\th#1 \enonce#2\endth{\global\def\currenvir{th}
   \medbreak\proclaim{#1}{\it #2}\medskip}

\long\def\comment#1\endcomment{}


\def\isinlabellist#1\of#2{\notfoundtrue%
   {\def\given{#1}%
    \def\\##1{\def\next{##1}%
    \lop\next\to\za\lop\next\to\zb%
    \ifx\za\given{\zb\global\notfoundfalse}\fi}#2}%
    \ifnotfound{\immediate\write16%
                 {Warning - [Page \the\pageno] {#1} No reference found}}%
                \fi}%
\def\ref#1{\ifx\labellist\empty{\immediate\write16
                 {Warning - No references found at all.}}
               \else{\isinlabellist{#1}\of\labellist}\fi}

\def\newlabel#1#2{\rightappenditem{\\{#1}\\{#2}}\to\labellist}
\def\labellist{}

\def\label#1{%
  \def\given{th}%
  \ifx\given\currenvir%
{\hwrite\lbl{\string\newlabel{#1}{\number\secno.\number\prmno}}}\fi%

\def\given{section}%
  \ifx\given\currenvir%
{\hwrite\lbl{\string\newlabel{#1}{\number\secno}}}\fi%

\def\given{subsection}%
  \ifx\given\currenvir%
{\hwrite\lbl{\string\newlabel{#1}{\number\secno.\number\prmno}}}\fi%

\def\given{formule}%
  \ifx\given\currenvir%
{\hwrite\lbl{\string\newlabel{#1}{\number\secno.\number\prmno}}}\fi%

\def\given{subsubsection}%
  \ifx\given\currenvir%
  {\hwrite\lbl{\string%
   
\newlabel{#1}{\number\secno.\number\subsecno.\number\subsubsecno}}}\fi
  \ignorespaces}

\newwrite\lbl

\def\openall{\openout\lbl=\jobname.lbl}
\def\closeall{\closeout\lbl}

\newread\testfile
\def\lookatfile#1{\openin\testfile=\jobname.#1
    \ifeof\testfile{\immediate\openout\nameuse{#1}\jobname.#1
                    \write\nameuse{#1}{}
                    \immediate\closeout\nameuse{#1}}\fi%
    \immediate\closein\testfile}%

\def\begin{\lookatfile{lbl}
           \input\jobname.lbl
           \openall}
\let\bye\end
\def\end{\closeall\bye}


\font\eightrm=cmr8         \font\eighti=cmmi8
\font\eightsy=cmsy8        \font\eightbf=cmbx8
\font\eighttt=cmtt8        \font\eightit=cmti8
\font\eightsl=cmsl8        \font\sixrm=cmr6
\font\sixi=cmmi6           \font\sixsy=cmsy6
\font\sixbf=cmbx6


\font\tengoth=eufm10       \font\tenbboard=msbm10
\font\eightgoth=eufm8      \font\eightbboard=msbm8
\font\sevengoth=eufm7      \font\sevenbboard=msbm7
\font\sixgoth=eufm6        \font\fivegoth=eufm5

\font\tencyr=wncyr10       
\font\eightcyr=wncyr8      
\font\sevencyr=wncyr7      
\font\sixcyr=wncyr6


\skewchar\eighti='177 \skewchar\sixi='177
\skewchar\eightsy='60 \skewchar\sixsy='60


\newfam\gothfam           \newfam\bboardfam
\newfam\cyrfam

\def\tenpoint{%
  \textfont0=\tenrm \scriptfont0=\sevenrm \scriptscriptfont0=\fiverm
  \def\rm{\fam\z@\tenrm}%
  \textfont1=\teni  \scriptfont1=\seveni  \scriptscriptfont1=\fivei
  \def\oldstyle{\fam\@ne\teni}\let\old=\oldstyle
  \textfont2=\tensy \scriptfont2=\sevensy \scriptscriptfont2=\fivesy
  \textfont\gothfam=\tengoth \scriptfont\gothfam=\sevengoth
  \scriptscriptfont\gothfam=\fivegoth
  \def\goth{\fam\gothfam\tengoth}%
  \textfont\bboardfam=\tenbboard \scriptfont\bboardfam=\sevenbboard
  \scriptscriptfont\bboardfam=\sevenbboard
  \def\bb{\fam\bboardfam\tenbboard}%
 \textfont\cyrfam=\tencyr \scriptfont\cyrfam=\sevencyr
  \scriptscriptfont\cyrfam=\sixcyr
  \def\cyr{\fam\cyrfam\tencyr}%
  \textfont\itfam=\tenit
  \def\it{\fam\itfam\tenit}%
  \textfont\slfam=\tensl
  \def\sl{\fam\slfam\tensl}%
  \textfont\bffam=\tenbf \scriptfont\bffam=\sevenbf
  \scriptscriptfont\bffam=\fivebf
  \def\bf{\fam\bffam\tenbf}%
  \textfont\ttfam=\tentt
  \def\tt{\fam\ttfam\tentt}%
  \abovedisplayskip=12pt plus 3pt minus 9pt
  \belowdisplayskip=\abovedisplayskip
  \abovedisplayshortskip=0pt plus 3pt
  \belowdisplayshortskip=4pt plus 3pt 
  \smallskipamount=3pt plus 1pt minus 1pt
  \medskipamount=6pt plus 2pt minus 2pt
  \bigskipamount=12pt plus 4pt minus 4pt
  \normalbaselineskip=12pt
  \setbox\strutbox=\hbox{\vrule height8.5pt depth3.5pt width0pt}%
  \let\bigf@nt=\tenrm       \let\smallf@nt=\sevenrm
  \normalbaselines\rm}

\def\eightpoint{%
  \textfont0=\eightrm \scriptfont0=\sixrm \scriptscriptfont0=\fiverm
  \def\rm{\fam\z@\eightrm}%
  \textfont1=\eighti  \scriptfont1=\sixi  \scriptscriptfont1=\fivei
  \def\oldstyle{\fam\@ne\eighti}\let\old=\oldstyle
  \textfont2=\eightsy \scriptfont2=\sixsy \scriptscriptfont2=\fivesy
  \textfont\gothfam=\eightgoth \scriptfont\gothfam=\sixgoth
  \scriptscriptfont\gothfam=\fivegoth
  \def\goth{\fam\gothfam\eightgoth}%
  \textfont\cyrfam=\eightcyr \scriptfont\cyrfam=\sixcyr
  \scriptscriptfont\cyrfam=\sixcyr
  \def\cyr{\fam\cyrfam\eightcyr}%
  \textfont\bboardfam=\eightbboard \scriptfont\bboardfam=\sevenbboard
  \scriptscriptfont\bboardfam=\sevenbboard
  \def\bb{\fam\bboardfam}%
  \textfont\itfam=\eightit
  \def\it{\fam\itfam\eightit}%
  \textfont\slfam=\eightsl
  \def\sl{\fam\slfam\eightsl}%
  \textfont\bffam=\eightbf \scriptfont\bffam=\sixbf
  \scriptscriptfont\bffam=\fivebf
  \def\bf{\fam\bffam\eightbf}%
  \textfont\ttfam=\eighttt
  \def\tt{\fam\ttfam\eighttt}%
  \abovedisplayskip=9pt plus 3pt minus 9pt
  \belowdisplayskip=\abovedisplayskip
  \abovedisplayshortskip=0pt plus 3pt
  \belowdisplayshortskip=3pt plus 3pt 
  \smallskipamount=2pt plus 1pt minus 1pt
  \medskipamount=4pt plus 2pt minus 1pt
  \bigskipamount=9pt plus 3pt minus 3pt
  \normalbaselineskip=9pt
  \setbox\strutbox=\hbox{\vrule height7pt depth2pt width0pt}%
  \let\bigf@nt=\eightrm     \let\smallf@nt=\sixrm
  \normalbaselines\rm}

\tenpoint


\def\pc#1{\bigf@nt#1\smallf@nt}         \def\pd#1 {{\pc#1} }

\catcode`\;=\active
\def;{\relax\ifhmode\ifdim\lastskip>\z@\unskip\fi
\kern\fontdimen2  -1.2 \fontdimen3 \string;}

\catcode`\:=\active
\def:{\relax\ifhmode\ifdim\lastskip>\z@\unskip\fi\penalty\@M\ \fi\string:}

\catcode`\!=\active
\def!{\relax\ifhmode\ifdim\lastskip>\z@
\unskip\fi\kern\fontdimen2  -1.1 \fontdimen3 \string!}

\catcode`\?=\active
\def?{\relax\ifhmode\ifdim\lastskip>\z@
\unskip\fi\kern\fontdimen2  -1.1 \fontdimen3 \string?}

\def\^#1{\if#1i{\accent"5E\i}\else{\accent"5E #1}\fi}
\def\"#1{\if#1i{\accent"7F\i}\else{\accent"7F #1}\fi}

\frenchspacing


\newtoks\auteurcourant      \auteurcourant={\hfil}
\newtoks\titrecourant       \titrecourant={\hfil}

\newtoks\hautpagetitre      \hautpagetitre={\hfil}
\newtoks\baspagetitre       \baspagetitre={\hfil}

\newtoks\hautpagegauche     
\hautpagegauche={\eightpoint\rlap{\folio}\hfil\the\auteurcourant\hfil}
\newtoks\hautpagedroite     
\hautpagedroite={\eightpoint\hfil\the\titrecourant\hfil\llap{\folio}}

\newtoks\baspagegauche      \baspagegauche={\hfil} 
\newtoks\baspagedroite      \baspagedroite={\hfil}

\newif\ifpagetitre          \pagetitretrue




\def\raggedbottom{\topskip 10pt plus 36pt\r@ggedbottomtrue}




\def\ctexte#1\endctexte{%
  \hbox{$\vcenter{\halign{\hfill##\hfill\crcr#1\crcr}}$}}


\long\def\ctitre#1\endctitre{%
    \ifdim\lastskip<24pt\vskip-\lastskip\bigbreak\bigbreak\fi
  		\vbox{\parindent=0pt\leftskip=0pt plus 1fill
          \rightskip=\leftskip
          \parfillskip=0pt\bf#1\par}
    \bigskip\nobreak}

\let\+=\tabalign

\def\signature#1\endsignature{\vskip 15mm minus 5mm\rightline{\vtop{#1}}}

\mathcode`A="7041 \mathcode`B="7042 \mathcode`C="7043 \mathcode`D="7044
\mathcode`E="7045 \mathcode`F="7046 \mathcode`G="7047 \mathcode`H="7048
\mathcode`I="7049 \mathcode`J="704A \mathcode`K="704B \mathcode`L="704C
\mathcode`M="704D \mathcode`N="704E \mathcode`O="704F \mathcode`P="7050
\mathcode`Q="7051 \mathcode`R="7052 \mathcode`S="7053 \mathcode`T="7054
\mathcode`U="7055 \mathcode`V="7056 \mathcode`W="7057 \mathcode`X="7058
\mathcode`Y="7059 \mathcode`Z="705A
 
\def\spacedmath#1{\def\packedmath##1${\bgroup\mathsurround=0pt ##1\egroup$}%
\mathsurround#1 \everymath={\packedmath}\everydisplay={\mathsurround=0pt }}
 
\def\nospacedmath{\mathsurround=0pt \everymath={}\everydisplay={} }


\def\decale#1{\smallbreak\hskip 28pt\llap{#1}\kern 5pt}
\def\decaledecale#1{\smallbreak\hskip 34pt\llap{#1}\kern 5pt}
\def\puce{\smallbreak\hskip 6pt{$\scriptstyle\bullet$}\kern 5pt}


\def\displaylinesno#1{\displ@y\halign{
\hbox to\displaywidth{$\@lign\hfil\displaystyle##\hfil$}&
\llap{$##$}\crcr#1\crcr}}


\def\ldisplaylinesno#1{\displ@y\halign{ 
\hbox to\displaywidth{$\@lign\hfil\displaystyle##\hfil$}&
\kern-\displaywidth\rlap{$##$}\tabskip\displaywidth\crcr#1\crcr}}


\def\eqalign#1{\null\,\vcenter{\openup\jot\m@th\ialign{
\strut\hfil$\displaystyle{##}$&$\displaystyle{{}##}$\hfil
&&\quad\strut\hfil$\displaystyle{##}$&$\displaystyle{{}##}$\hfil
\crcr#1\crcr}}\,}


\def\system#1{\left\{\null\,\vcenter{\openup1\jot\m@th
\ialign{\strut$##$&\hfil$##$&$##$\hfil&&
        \enskip$##$\enskip&\hfil$##$&$##$\hfil\crcr#1\crcr}}\right.}


\let\@ldmessage=\message

\def\message#1{{\def\pc{\string\pc\space}%
                \def\'{\string'}\def\`{\string`}%
                \def\^{\string^}\def\"{\string"}%
                \@ldmessage{#1}}}

\def\diagram#1{\def\normalbaselines{\baselineskip=0pt
\lineskip=10pt\lineskiplimit=1pt}   \matrix{#1}}


\def\up#1{\raise 1ex\hbox{\smallf@nt#1}}


\def\cf{{\it cf}} 

\def\qed{\raise -2pt\hbox{\vrule\vbox to 10pt{\hrule width 4pt
                 \vfill\hrule}\vrule}}

\def\cqfd{\unskip\penalty 500\quad\vrule height 4pt depth 0pt width 4pt\medbreak}

\def\virg{\raise .4ex\hbox{,}}   


\def\build#1_#2^#3{\mathrel{
\mathop{\kern 0pt#1}\limits_{#2}^{#3}}}


\def\boxit#1#2{%
\setbox1=\hbox{\kern#1{#2}\kern#1}%
\dimen1=\ht1 \advance\dimen1 by #1 \dimen2=\dp1 \advance\dimen2 by #1 
\setbox1=\hbox{\vrule height\dimen1 depth\dimen2\box1\vrule}%
\setbox1=\vbox{\hrule\box1\hrule}%
\advance\dimen1 by .4pt \ht1=\dimen1 
\advance\dimen2 by .4pt \dp1=\dimen2  \box1\relax}

\def\date{\the\day\ \ifcase\month\or janvier\or f\'evrier\or mars\or
avril\or mai\or juin\or juillet\or ao\^ut\or septembre\or octobre\or
novembre\or d\'ecembre\fi \ {\old \the\year}}

\def\dateam{\ifcase\month\or January\or February\or March\or
April\or May\or June\or July\or August\or September\or October\or
November\or December\fi \ \the\day ,\ \the\year}

\def\moins{\mathop{\hbox{\vrule height 3pt depth -2pt
width 5pt}}}
\def\crog{{\vrule height 2.57mm depth 0.85mm width 0.3mm}\kern -0.36mm
[}

\def\crod{]\kern -0.4mm{\vrule height 2.57mm depth 0.85mm
width 0.3 mm}}
\def\moins{\mathop{\hbox{\vrule height 3pt depth -2pt
width 5pt}}}

\def\rond{\kern 1pt{\scriptstyle\circ}\kern 1pt}

\def\diagram#1{\def\normalbaselines{\baselineskip=0pt
\lineskip=10pt\lineskiplimit=1pt}   \matrix{#1}}

\def\hfl#1#2{\nospacedmath\smash{\mathop{\hbox to
12mm{\rightarrowfill}}\limits^{\scriptstyle#1}_{\scriptstyle#2}}}

\def\ghfl#1#2{\nospacedmath\smash{\mathop{\hbox to
25mm{\rightarrowfill}}\limits^{\scriptstyle#1}_{\scriptstyle#2}}}

\def\phfl#1#2{\nospacedmath\smash{\mathop{\hbox to
8mm{\rightarrowfill}}\limits^{\scriptstyle#1}_{\scriptstyle#2}}}

\def\pvfl#1#2{\llap{$\scriptstyle#1$}\left\downarrow\vbox to 4mm{}\right.\rlap{$\scriptstyle#2$}}

\def\va{vari\'et\'e ab\'elienne}

\def\pa{\S\kern.15em}

\def\Z{\hbox{\bf Z}}
\def\N{\hbox{\bf N}}
\def\P{\hbox{\bf P}}

\def\Q{\hbox{\bf Q}}
\def\C{\hbox{\bf C}}

\def\cad{c'est-\`a-dire}

\def\Card{\mathop{\rm Card}\nolimits}

\def\cf{{\it cf.\/}}

\def\codim{\mathop{\rm codim}\nolimits}

\def\dra{\ra\kern -3mm\ra}

\def\gr{\mathop{\rm gr}\nolimits}

\def\indp{\par\hskip0.5cm}
\def\isom{\simeq}

\def\Ker{\mathop{\rm Ker}\nolimits}
\def\ldra{\lra\kern -3mm\ra}

\def\loc{{\it loc.cit.\/}}
\def\long{\mathop{\rm long}\nolimits}
\def\lra{\longrightarrow}
\def\llra{\nospacedmath\hbox to 10mm{\rightarrowfill}}
\def\lllra{\nospacedmath\hbox to 15mm{\rightarrowfill}}

\def\Pic{\mathop{\rm Pic}\nolimits}

\def\ra{\rightarrow}

\def\Sym{\mathop{\rm Sym}\nolimits}

\def\theo{th\'eor\`eme}

\def\tv{\tvi\vrule}
\def\tvi{\vrule height 12pt depth 5pt width 0pt}
\def\tx{\kern -1.5pt -}

\def\cc#1{\hfill\kern .7em#1\kern .7em\hfill}

\def\note#1#2{\footnote{\parindent
0.4cm$^#1$}{\vtop{\eightpoint\baselineskip12pt\hsize15.5truecm
\noindent #2}}\parindent 0cm}

\def\og{\leavevmode\raise.3ex\hbox{$\scriptscriptstyle\langle\!\langle$}}
\def\fg{\leavevmode\raise.3ex\hbox{$\scriptscriptstyle\,\rangle\!\rangle$}}

\catcode`\@=12

\showboxbreadth=-1  \showboxdepth=-1


\spacedmath{1.5pt}\parskip=1mm
\baselineskip=14pt
\font\eightrm=cmr10 at 8pt
\overfullrule=0pt
\input amssym.def 
\input amssym
\def\dra{\dashrightarrow} 

\def\eps{\epsilon}
\def\phi{\varphi}
\def\ot{\otimes}
\def\otc{{\mathord{\otimes\cdots\otimes }}}

\def\we{\bigwedge}

\def\a{{\alpha}}
\def\b{{\beta}}
\def\d{{\delta}}
\def\e{{\varepsilon}}
\def\k{{\kappa}}
\def\l{{\lambda}}

\def\s{{\sigma}}
\def\t{{\tau}}

\def\cS{{\cal S}}

\def\cI{{\cal I}}

\def\cO{{\cal O}}

\def\algc{alg\'ebriquement clos}
\def\lin{lin\'eaire}

\begin
\null\bigskip
\centerline {\bf SCH\'EMAS DE FANO} 
 \medskip
\centerline{{\bf Olivier Debarre\note{1}{\baselineskip=3truemm\rm Financ\'e en
partie par le Projet Europ\'een HCM \og Al\-ge\-braic
Geometry in Europe\fg\ (AGE), Contrat CHRXCT-940557.} et Laurent Manivel}}

\vskip 1cm

\ind Soient $k$ un corps \algc\ et $X$ un sous-sch\'ema d'un espace projectif $\P^n_k$; on
appelle sch\'ema de Fano, et l'on note $F_r(X)$, le sous-sch\'ema de la  grassmannienne
$G(r,\P^n_k)$ qui param\`etre les espaces
\lin s de dimension $r$ contenus dans $X$. Ces sch\'emas ont une longue histoire ([F], [vW], [AK], 
[BVV], [B1], [PS], [K],     [ELV], [BV]) mais il ne semble pas exister dans la
litt\'erature d'\'enonc\'e g\'en\'eral sur leurs propri\'et\'es, m\^eme les plus simples comme la
connexit\'e, valable en toute caract\'eristique. Un des buts de cet
article est de rassembler sous une r\'ef\'erence commune des faits g\'en\'eraux sur ces
sch\'emas.  

\ind Apr\`es un paragraphe de notations, on obtient dans le \S 2, en
se basant sur les id\'ees de [K],  notre premier r\'esultat
principal: pour une intersection compl\`ete g\'en\'erale $X$, le  sch\'ema de Fano $F_r(X)$
{\it est non vide et lisse de la dimension attendue
$\d$ lorsque celle-ci est positive, et connexe lorsque} $\d>0$. Dans le
\S 3, on applique des r\'esultats de [D] et [S] pour calculer  certains groupes d'homotopie de
$F_r(X)$. Par ailleurs, le sch\'ema
$F_r(X)$ est le lieu des z\'eros d'une section d'un fibr\'e vectoriel sur la grassmannienne
$G(r,\P^n)$; lorsqu'il a la dimension attendue $\delta$, son id\'eal admet une r\'esolution par un
complexe de Koszul. Un \theo\ d'annulation pour certains fibr\'es vectoriels sur $G(r,\P^n)$ (prop.
\ref{annul}) nous permet de montrer notre second r\'esultat principal, \`a savoir {\it un \theo\ de
type Lefschetz, qui permet d'obtenir, pour
 $k=\C$,   les nombres de Hodge  $h^{p,q}(F_r(X))$  pour
$p+q$ assez petit} (inf\'erieur \`a $\dim X-2r$ pour $n$ grand). Apr\`es avoir r\'edig\'e cette
partie, nous nous sommes rendus compte que Borcea avait d\'ej\`a utilis\'e le \theo\ d'annulation
de Bott dans ce cadre (il obtient entre autres les r\'esultats du \S 2 en caract\'eristique nulle).

\ind Les m\^emes m\'ethodes permettent d'\'etudier dans le \S 4 la restriction\break
$H^0(G(r,\P^n),\cO(l))\lra H^0(F_r(X),\cO(l))$; on montre que pour $n$ assez grand, $F_r(X)$
est projectivement normal dans $G(r,\P^n)$, et que toute \'equation de $F_r(X)$ est de 
degr\'e au moins \'egal \`a une \'equation de $X$ dans $\P^n$. On donne aussi une formule
explicite pour le calcul du degr\'e des sch\'emas $F_r(X)$: c'est le coefficient d'un mon\^ome
particulier dans un polyn\^ome explicite en
$r+1$ variables. On donne quelques exemples de ce calcul pour des hypersurfaces de bas degr\'e.

\ind On s'int\'eresse ensuite aux sous-sch\'emas de $F_r(X)$ qui param\`etrent les $r$\tx plans
contenant un $r_0$\tx plan fix\'e; le \theo\ principal du \S 5 g\'en\'eralise les r\'esultats
analogues du \S 2 dans ce cadre. On en d\'eduit que $F_r(X)$ est s\'eparablement unir\'egl\'e en
droites pour $n$ assez grand, ce qui nous permet dans le \S 6 d'adapter des id\'ees de [K] pour
montrer, toujours pour $n$ assez grand, que le groupe  de Chow rationnel  des $1$\tx cycles sur un
sch\'ema  $F_r(X)$ est de rang
$1$. Il est tentant de g\'en\'eraliser une conjecture de  Srinivas et Paranjape
([P]) de la fa\c con suivante: pour
$n$ assez grand, les groupes de Chow rationnels de basse dimension de $F_r(X)$ devraient \^etre ceux
de la grassmannienne ambiante $G(r,\P^n)$.

\ind Dans le \S 7, on d\'emontre, comme conjectur\'e dans
[BVV], que {\it les sch\'emas de Fano g\'en\'eriques sont unirationnels pour $n$ assez grand}. On se
ram\`ene pour cela
\`a un r\'esultat de Predonzan ([Pr]), pr\'ecis\'e dans l'article  [PS], qui fournit un crit\`ere
explicite pour l'unirationa\-li\-t\'e d'une intersection compl\`ete dans un espace projectif.
Les bornes obtenues sont explicites, mais tr\`es grandes; par exemple, on  montre que la vari\'et\'e
des droites contenues dans une hypersurface cubique de
$\P^n$ est unirationnelle pour $n\ge 433$ (alors que c'est d\'ej\`a une vari\'et\'e de Fano pour
$n\ge 6$).

\section{Notations}

\ind Soient $k$ un corps alg\'ebriquement clos et $V$ un $k$\tx espace vectoriel de dimension
$n+1$.   Pour toute suite finie ${\bf
d}=(d_1,\ldots,d_s)$ d'entiers positifs, et tout entier positif $r$, on note $|{\bf
d}|=\sum_{i=1}^sd_i$, puis ${\bf
d}+r=(d_1+r,\ldots,d_s+r)$ et ${{\bf d}\choose r}={\sum_{i=1}^s}{d_i\choose
r} $. On pose  $\Sym ^{\bf
d}V^*=\bigoplus_{i=1}^s\Sym ^{d_i}V^*$, espace vectoriel que l'on notera aussi $\Gamma_{{\bf P}
V}({\bf d})$. Enfin, si
${\bf f}=(f_1\ldots,f_s)$ est un \'el\'ement non nul de
$\Sym ^{\bf d}V^*$, on note $X_{\bf f}$ le sous-sch\'ema de $\P V$ d'\'equations $f_1=\cdots=f_s=0$; on
dira d'un tel sch\'ema qu'il est {\it d\'efini par des \'equations de degr\'e} ${\bf d}$.

\ind  On pose ensuite 
$$ \delta(n,{\bf d},r)=(r+1)(n-r)- {{\bf d}+r\choose r}$$
et $\delta_-(n,{\bf d},r)=\min\{ \delta (n,{\bf d},r), n-2r-s\}$, que l'on \'ecrira simplement
$\delta$ et $\delta_-$ lorsque qu'aucune confusion ne sera \`a craindre. 

 \section{Dimension, lissit\'e et connexit\'e}

\ind On montre dans ce num\'ero que les sch\'emas de Fano d'un sous-sch\'ema $X$
de $\P^n_k$
d\'efini par des \'equations de degr\'es ${\bf d}=(d_1,\ldots,d_s)$ sont lisses de la dimension
attendue pour $X$ g\'en\'erale, et connexes lorsque cette dimension est strictement positive. Divers
cas particuliers du \theo\ suivant \'etaient d\'ej\`a connus: citons par exemple [BVV], qui traite
le cas $k=\C$ et $r=s=1$; [P], [Mu] et [PS], qui d\'emontrent b); [B1], qui d\'emontre le \theo\
lorsque $k$ est de caract\'eristique nulle; et [K], qui traite le cas
$r=s=1$ (th. 4.3, p. 266), et dont nous empruntons les id\'ees. Lorsque $k=\C$, une d\'emonstration
compl\`etement diff\'erente d\'ecoule de celle du th\'eor\`eme \ref{leff} (\cf\ 
rem. \ref{rem}.1).

\ind Pour appliquer le \theo , il est utile de noter que {\it lorsque ${\bf d}\ne (2)$, l'entier
  $\delta(n,{\bf d},r)$ est positif (resp. strictement positif)  si et seulement si
$\delta_-(n,{\bf d},r)$ l'est.}   

\th
\label{fano}
Th\'eor\`eme
\enonce
Soient $X$ un sous-sch\'ema de
$\P^n_k$ d\'efini par des \'equations de degr\'e ${\bf d}$,  et $F_r(X)$ le sch\'ema de
Fano des
$r$\tx plans contenus dans $X$.

\indp {\rm a)} Lorsque $\delta_-(n,{\bf d},r)< 0$, le sch\'ema $F_r(X)$ est vide pour $X$
g\'en\'erale.

\indp {\rm b)} Lorsque $\delta_-(n,{\bf d},r)\ge 0$, le sch\'ema $F_r(X)$ est non vide; il est
lisse de dimension
$\delta(n,{\bf d},r)$ pour $X$ g\'en\'erale.

\indp {\rm c)} Lorsque $\delta_-(n,{\bf d},r)>0$, le sch\'ema $F_r(X)$ est connexe.
\endth

\ind Consid\'erons la
vari\'et\'e d'incidence 
$$I_r=\{ ([{\bf f}],\Lambda)\in
 \P \Sym ^{\bf d}V^*\times G(r,\P^n) \mid \Lambda\i X_{\bf f}\}\ ,$$ et les 
projections
$p_r:I_r\ra
 \P \Sym ^{\bf d}V^*$ (dont la fibre au-dessus de $[{\bf f}]$ s'identifie \`a $F_r(X_{\bf f})$) et
$q:I_r\ra G(r,\P^n)$. Etant donn\'e un $r$\tx plan $\Lambda=\P W$, la fibre 
$q^{-1}([\Lambda])$ est l'espace projectif associ\'e au noyau du morphisme
surjectif
$\Sym ^{\bf d}V^*\ra \Sym ^{\bf d}W^*$. Elle est donc de codimension
${{\bf d}+r\choose r}$ dans
$\P \Sym ^{\bf d}V^*$, de sorte que $I_r$ est irr\'eductible lisse de codimension 
${{\bf d}+r\choose r}$ dans $
\P \Sym ^{\bf d}V^*\times G(r,\P^n)$.

\ind On note $Z_r$ le ferm\'e des points de $I_r$ o\`u $p_r$ n'est pas lisse, et $\Delta_r$ l'image
de $Z_r$ par $p_r$ (avec la convention $\Delta_{-1}=\varnothing$). Soit $\Lambda$ un $r$\tx plan,
d'\'equations $x_{r+1}=\cdots=x_n=0$ dans $\P V$; pour tout entier $m\ge 0$, on note
${\cal B}_m$ la base $\{ {\bf x}^J\mid J\i
\{0,\ldots,r\}\ ,\ \Card(J)=m\}$ de l'espace vectoriel $\Gamma_\Lambda(m)$; on note aussi
${\cal B}_{\bf d}$ la base
$\iota_1({\cal B}_{d_1})\cup\cdots\cup\iota_s({\cal B}_{d_s})$ de l'espace vectoriel $\Gamma_\Lambda({\bf
d})$ (o\`u $\iota_j$ est l'injection canonique de
$\Gamma_\Lambda(d_j)$ dans $\Gamma_\Lambda({\bf d})$). 

\th
\label{lisse}
Lemme
\enonce Pour qu'un point
$([{\bf f}],[\Lambda])$ de
$I_r$ soit dans $Z_r$, il faut et il suffit que le morphisme
$\alpha:\Gamma_\Lambda(1)^{n-r}\ra \Gamma_\Lambda({\bf d})$ d\'efini par 
$$\alpha(h_{r+1},\ldots,h_n)=\Bigl( \sum_{i=r+1}^nh_i \Bigl({\partial
f_1\over\partial x_i}\Bigr) _{\displaystyle\vert_{\scriptstyle\Lambda}},\ldots,
\sum_{i=r+1}^nh_i \Bigl({\partial
f_s\over\partial x_i}\Bigr) _{\displaystyle\vert_{\scriptstyle\Lambda}}\Bigr)$$
ne soit pas surjectif.
\endth

\ind Cela r\'esulte d'un calcul explicite fait dans [BVV] dans le cas $r=s=1$.\cqfd

\ind {\it Lorsque $X$ est lisse de codimension $r$ le long de}
$\Lambda$, on a une suite exacte 
$$0\lra N_{\Lambda/X}\lra {\cal
O}_\Lambda(1)^{n-r}\buildrel{u}\over{\lra} \bigoplus_{i=1}^s{\cal O}_\Lambda(d_i)\lra
0
\ ,$$  et le morphisme $\alpha$ n'est
autre que $H^0(u)$ (\cf\ [BVV], prop. 3 et [K], p. 267 dans le cas $r=s=1$). La condition du
lemme est donc \'equivalente dans ce cas \`a l'annulation de
$H^1(\Lambda,N_{\Lambda/X})$.

\ind Soit
$\mu: \Gamma_\Lambda(1)\times \Gamma_\Lambda({\bf d}-1)\ra \Gamma_\Lambda({\bf d})$
le morphisme de multiplication, d\'efini par\break $\mu(h,g_1,\ldots,g_s)=(hg_1,\ldots,hg_s)$. Si $H$
est un hyperplan de $\Gamma_\Lambda({\bf d})$, on note 
$\mu^{-1}(H)$ l'ensemble $\{\ g\in \Gamma_\Lambda({\bf d}-1)\mid
\mu(\Gamma_\Lambda(1)\times\{ g\} )\i
H\
\}$. 

\ind On peut r\'e\'enoncer le lemme \ref{lisse} de la fa\c con suivante: soit ${\cal Z}$ le
sous-ensemble de $q^{-1}([\Lambda])\times\P  \Gamma_\Lambda({\bf d})$ form\'e des couples
$([{\bf f}],[\ell])$ tels que
$$\Bigl(  \Bigl({\partial
f_1\over\partial x_i}\Bigr) _{\displaystyle\vert_{\scriptstyle\Lambda}},\ldots,
 \Bigl({\partial
f_s\over\partial x_i}\Bigr) _{\displaystyle\vert_{\scriptstyle\Lambda}}\Bigr)$$
soit dans $\mu^{-1}(\Ker(\ell))$ pour tout $i=r+1,\ldots,n$; alors $Z_r\cap
q^{-1}([\Lambda])$ s'identifie \`a la premi\`ere projection de ${\cal Z}$. Pour tout entier
$h$, notons ${\cal L}_h$ l'ensemble
 des formes lin\'eaires $\ell$ sur $\Gamma_\Lambda({\bf d})$
 v\'erifiant
$\codim \mu^{-1}(\Ker(\ell)) =h$, et ${\cal Z}_h$ l'ensemble des \'el\'ements 
$([{\bf f}],[\ell])$ de ${\cal Z}$ avec $\ell\in {\cal L}_h$. On peut \'ecrire
$f_i=\sum_{j=r+1}^nx_jf_{ij}$, avec 
$f_{ij}\vert_\Lambda=\bigl({\partial
f_i\over\partial x_j}\bigr) _{\displaystyle\vert_{\scriptstyle\Lambda}}$, de sorte que 
$$\codim_{q^{-1}([\Lambda])}pr_1({\cal Z}_h)\ge h(n-r)-\dim\P {\cal
L}_h\leqno{\formule}\hbox{\label{minor1}}$$ et
$$\codim_{I_r}Z_r=\codim_{q^{-1}([\Lambda])}pr_1({\cal Z})\ge\min_{1\le h\le r+1}[
h(n-r)-\dim\P {\cal L}_h]\ .\leqno{\formule}\hbox{\label{minor2}}$$

\subsection Soit $\ell$ une forme lin\'eaire sur $\Gamma_\Lambda({\bf d})$. Soit $M$ la matrice
 \`a coefficients dans $\Gamma_\Lambda({\bf d})$ de la forme bilin\'eaire $\mu$ dans les bases
${\cal B}_1$ et ${\cal B}={\cal B}_{{\bf d}-1}$. Pour qu'un \'el\'ement $g=\sum_{b\in{\cal B}} g_bb$ de
$\Gamma_\Lambda({\bf d}-1)$ soit dans $\mu^{-1}(\Ker(\ell))$, il faut et il suffit que $\sum_b
g_b\ell(x_ib)$ soit nul pour tout $i=0,\ldots,r$, de sorte que {\it la
codimension de $\mu^{-1}(\Ker(\ell))$ dans $\Gamma_\Lambda({\bf d}-1)$ est le rang de la matrice} 
$\ell(M)$.\label{forme}

\th
\label{propre}
Lemme
\enonce
Soient $([{\bf f}],[\Lambda])$ un \'el\'ement de $Z_r\moins
p_r^{-1}(\Delta_{r-1})$ et $\ell$ une forme lin\'eaire non nulle sur $\Gamma_\Lambda({\bf d})$, qui s'annule
sur l'image de $\alpha$. Alors
$\mu^{-1}(\Ker(\ell))$ est de codimension $r+1$ dans $\Gamma_\Lambda({\bf d}-1)$.
\endth

\ind Proc\'edons par l'absurde en supposant que la matrice $\ell(M)$ d\'efinie ci-dessus ne soit
pas de rang maximal. Quitte \`a effectuer un changement lin\'eaire de coordonn\'ees, on peut
supposer $\ell(x_rb)=0$ pour tout $b$ dans ${\cal B}$, de sorte que si $\Lambda'$ est l'hyperplan de
$\Lambda$ d\'efini par $x_r=0$, la forme lin\'eaire $\ell$ provient d'une forme lin\'eaire $\ell'$
sur $ \Gamma_{\Lambda'}({\bf d})$. Si $\alpha':\Gamma_{\Lambda'}(1)^{n-r+1}\ra
\Gamma_{\Lambda'}({\bf d})$ est le morphisme associ\'e au point $([{\bf f}],[\Lambda'])$ de
$I_{r-1}$ d\'efini dans le lemme
\ref{lisse}, $\ell'$ s'annule sur
$\alpha'(\{ 0\}\times \Gamma_{\Lambda'}(1)^{n-r})$. Comme la restriction de ${\partial f_i\over
\partial x_r}$ \`a $\Lambda'$ est nulle pour tout $i$, la forme lin\'eaire $\ell'$ s'annule sur
toute l'image de $\alpha'$, ce qui contredit  l'hypoth\`ese
$[{\bf f}]\notin \Delta_{r-1}$.\cqfd

\ind En d'autres termes, $q^{-1}([\Lambda])\cap\bigl( Z_r\moins
p_r^{-1}(\Delta_{r-1})\bigr)$ est contenu dans $pr_1({\cal Z}_{r+1})$, et (\ref{minor1})
entra\^ine
$$\nospacedmath\displaylines{\dim
(\overline{Z_r\moins
p_r^{-1}(\Delta_{r-1})})=\dim I_r-\codim_{I_r}
(\overline{Z_r\moins
p_r^{-1}(\Delta_{r-1})})\cr
\le\dim \P \Sym ^dV^*+\dim G(r,\P^n)-{{\bf d}+r\choose
r}-(r+1)(n-r)+\dim\P\Gamma_\Lambda({\bf d})
<\dim \P \Sym ^dV^*\ .\cr}
$$
\subsection Il en r\'esulte $\Delta_r\moins\Delta_{r-1}\ne \P \Sym ^dV^*$, d'o\`u
$\Delta_r\ne\P \Sym ^dV^*$ par r\'ecurrence sur $r$.\label{ferme}

\ind  On remarquera que nous avons en fait d\'emontr\'e que $\Delta_r$ a au plus une
composante irr\'eductible de plus que $\Delta_{r-1}$, \cad\ au plus $r+1$ composantes
irr\'eductibles. 

\th
\label{major}
Lemme
\enonce
Pour  $1\le h\le r+1$, la dimension de
${\cal L}_h$ est au plus $h(r-h+1)+{{\bf d}+h-1\choose h-1}$.
\endth

\ind On garde les notations de (\ref{forme}). Supposons les $h$ premi\`eres lignes de la
matrice $\ell(M)$
 lin\'eairement ind\'ependantes; on peut \'ecrire
$\ell(x_jb)=\sum_{i=0}^{h-1}a_{ij}\ell(x_ib)$, pour tous
$j=h,\ldots,r$ et
$b\in{\cal B}$, de sorte que les $\ell(b_i)$, pour $b_i={\bf x}^I$ dans ${\cal B}_{d_i}$,
peuvent s'exprimer en fonction de ceux pour lesquels $I\i \{0,\ldots,h-1\}$, et des
$h(r-h+1)$ coefficients $a_{ij}$.\cqfd

\ind L'in\'egalit\'e (\ref{minor2}) donne 
$$\codim_{I_r}Z_r\ge \min_{1\le h\le
r+1}\bigl[ h(n-2r+h-1)-{{\bf d}+h-1\choose h-1}\bigr]+1\ .$$
 \subsection Lorsque ${\bf d}\ne (2)$, on v\'erifie que l'expression entre crochets est une
fonction {\it concave} $\phi$ de $h$ sur $[1,+\infty [$; lorsque ${\bf d}=(2)$ et  $\delta_-\ge
0$, c'est une fonction {\it croissante}. On a dans chacun de ces cas\label{concave}\vskip-5mm 
$$\codim_{I_r}Z_r\ge \min \{ \phi(1),\phi(r+1)\}+1=\delta_-+1\
.$$

\ind Supposons $\delta_-<0$. Si ${\bf d}=(2)$, cela signifie $2r\ge n$; si une quadrique $X$
contient un $r$\tx plan $\Lambda$, \'ecrivons en gardant les m\^emes notations
$f=x_{r+1}\ell_{r+1}+\cdots+x_n\ell_n$, o\`u les $\ell_i$ sont des formes
\lin s. Comme $n-r\le r$, celles-ci ont un z\'ero commun sur $\Lambda$, qui est un point singulier de
$X$, ce qui ne peut se produire pour $X$ g\'en\'erale. Lorsque  ${\bf d}\ne (2)$, on a
$\delta<0$, d'o\`u
$\dim I_r<\dim \P \Sym ^dV^*$, et
$p_r$ n'est pas surjective; ceci   montre a) dans tous les cas. 

\ind Supposons $\delta_-\ge 0$; il existe d'apr\`es (\ref{concave}) un point
de $I_r$ en lequel $p_r$ est lisse. Cela entra\^ine que
$p_r$ est surjective, et que $F_r(X)$ est de dimension $\delta$ pour $X$
g\'en\'erale. Par (\ref{ferme}), $p_r$ est lisse au-dessus d'un ouvert dense de
$\P \Sym ^dV^*$, ce qui montre b).

\ind Supposons maintenant $\delta_->0$, et consid\'erons comme dans [BVV] la
factorisation de Stein $p_r:I_r\lra S\buildrel{\pi}\over{\lra}\P \Sym ^dV^*$ du
morphisme propre $p_r$. Si le morphisme $\pi$ est ramifi\'e, le \theo\ de
puret\'e entra\^ine que $Z_r$ contient l'image inverse d'un diviseur de $S$, ce qui
contredit l'estimation de (\ref{concave}). Il s'ensuit que
$\pi$ est \'etale, donc que c'est un isomorphisme puisque $\P \Sym ^dV^*$ est
simplement connexe. La vari\'et\'e $F_r(X)$ est donc connexe pour toute
hypersurface $X$, ce qui montre c).\cqfd

\ex{Remarques} 1) Soit $S$ le sous-fibr\'e tautologique sur $G(r,{\bf P}V)$. Tout \'el\'ement ${\bf
f}$ de  $\Sym ^{\bf d}V^*$ induit une section du fibr\'e  $\Sym ^dS$, dont
le lieu des z\'eros est le sch\'ema $F_r(X_{\bf f})$. La partie b) du \theo\ montre que lorsque 
$\delta_-(n,{\bf
d},r)\ge 0$, la classe de Chern $c_{\max}(\Sym ^{\bf d}S^*)$ est non nulle. On verra dans le \S 4 
comment expliciter cette classe de Chern dans l'anneau
de Chow de la grassmannienne. On remarque que lorsque ${\bf d}=(2)$ et que
$\delta_-<0\le\d$, le rang de $\Sym ^2S$ est plus petit que la dimension de $G(r,{\bf
P}^n)$, mais sa classe de Chern d'ordre maximal $2^{r+1}\sigma_{r+1,r,\ldots,1}$ est nulle (\cf\
[Fu], ex. 14.7.15).\label{chern}

 2) Toute quadrique lisse $X$ dans $\P^n$ est projectivement \'equivalente \`a la quadrique
d'\'equation $x_0x_1+x_2x_3+\cdots+x_{n-1}x_n=0$ si $n$ est impair, \`a la quadrique d'\'equation
$x_0x_1+x_2x_3+\cdots+x_{n-2}x_{n-1}+x_n^2=0$ si $n$ est pair. Le sch\'ema $F_r(X)$ est donc lisse
connexe d\`es que $\delta_->0$, \cad\ $n>2r+1$; on sait qu'il a deux composantes connexes si
$n=2r+1$.

\section{Groupes d'homotopie, groupes de cohomologie et groupe de Picard}

\ind Les r\'esultats de [D] et [S] permettent de calculer les groupes d'homotopie des sch\'emas
de Fano pour $n$ assez grand.

\th\label{pij} 
Proposition
\enonce
Soit $X$ un sous-sch\'ema de $\P^n_k$
d\'efini par des \'equations de degr\'e ${\bf d}$. On suppose $F_r(X)$ irr\'eductible de dimension
$\delta$.

\ind {\rm a)} Si $n\ge {2\over
r+1}{{\bf d}+r\choose r}+r+1$, le sch\'ema $F_r(X)$ est alg\'ebriquement simplement connexe,
topologiquement simplement connexe lorsque $k=\C$.

\ind {\rm b)} Lorsque $k=\C$ et que $F_r(X)$  est lisse, on a
$\pi_j\bigl( G(r,{\bf P}^n),F_r(X)\bigr)=0$ pour $n\ge 2{{\bf
d}+r\choose r}+j-1$. En particulier, si $n\ge 2{{\bf
d}+r\choose r}+2 $, le groupe de Picard de $F_r(X)$ est isomorphe \`a $\Z$, engendr\'e par la
classe de $\cO(1)$.
 \endth

\ind Le point b) est cons\'equence directe de [S]. Pour a), il suffit par [D], cor. 7.4
de montrer que  $[F_r(X)]\cdot[F_r(X)]\cdot[G(r,\P^{n-1})]$ est non nul dans $A(G(r,\P^n))$. Par
la remarque \ref{chern}, cette intersection est la classe de Chern de degr\'e maximal de $\Sym ^{\bf
d}S^*\oplus  \Sym ^{\bf d}S^*\oplus S^*$, et celle-ci est non nulle d\`es que $\delta(n-1,({\bf
d},{\bf d}),r)$ est positif, condition qui d\'ecoule de l'hypoth\`ese.\cqfd 

\ex{Remarques} 1)  On rappelle que $\pi_j(G(r,{\bf P}^n))\isom \pi_{j-1}(U(r+1))$ pour $j\le
2(n-r)$ ([H], chap. 7); si l'on suppose aussi  
$j\le 2(r+1)$, le th\'eor\`eme de p\'eriodicit\'e de Bott 
implique donc $\pi_j(G(r,{\bf P}^n))={\bf Z}$ ou $0$ selon que $j$ 
est pair ou impair. En g\'en\'eral, il peut cependant  
appara\^{\i}tre de la torsion (par exemple, $\pi_{11}(G(3,\P^n))=
{\bf Z}_2\oplus {\bf Z}_{120}$  si $n\ge 9$).\label{rempi}  

2) La remarque \ref{chern} montre que lorsque $F_r(X)$ est de dimension
$\delta$, on a 
$$\omega_{F_r(X)}\isom \omega_{G(r,{\bf P}^n)}\otimes \we^{\max}\Sym ^{\bf d}S^*\vert_{F_r(X)}
\isom {\cal O}_{F_r(X)}({{\bf d}+r\choose
r+1}-n-1)\ .$$
\ind En particulier, $F_r(X)$ est une vari\'et\'e de Fano lorsque $n\ge {{\bf d}+r\choose
r+1}$, donc simplement connexe lorsque $k=\C$ ([C1], [KMM1]). Cette borne est n\'eanmoins  
  moins bonne que celle de la prop. \ref{pij}.a) d\`es que l'un des $d_i$ est $\ge 3$..

\medskip
\ex{Exemple} Soit $X$ une hypersurface cubique lisse dans $\P^n_k$; par [BVV], prop. 5,
$F_1(X)$ est une vari\'et\'e {\it  lisse connexe} de dimension $2n-6$. La proposition
entra\^ine que $F_1(X)$ est simplement connexe pour $n\ge 6$. Lorsque $k=\C$, cela reste
vrai pour
$n=5$ ([BD], prop. 3), mais pas pour $n=4$ puisque $h^1(F_1(X),{\cal O}_{F_1(X)})=5$ ([AK], prop.
1.15).\label{BD}

\bigskip
\ind Passons maintenant au r\'esultat principal de ce num\'ero. On a vu en \ref{chern} que
$F_r(X)$ est le lieu des z\'eros d'une section d'un fibr\'e vectoriel sur la grassmannienne;
lorsqu'il a la codimension attendue, son id\'eal admet une r\'esolution par un
complexe de Koszul. Lorsque $k=\C$, on montre \`a l'aide du th\'eor\`eme de Borel-Weil-Bott ([Bo],
[De]) et des r\'esultats de [Ma1] et [Ma2] un \theo\ d'annulation (prop. \ref{annul}) qui nous
permettra
 de d\'eterminer certains
groupes de cohomologie des sch\'emas de Fano.

\th 
\label{leff}
Th\'eor\`eme
\enonce
Soit $X$ un sous-sch\'ema de $\P^n_{\bf C}$
d\'efini par des \'equations de degr\'e ${\bf d}$, tel que
 $F_r(X)$  soit lisse de dimension $\delta(n,{\bf d},r)$.  Le morphisme de restriction\break
$H^i(G(r,\P^n),\Q)\ra H^i(F_r(X),\Q)$
  est bijectif pour $i<\delta_-(n,{\bf d},r)$, 
 injectif pour
  $i=\delta_-(n,{\bf d},r)$. \endth

\ind En particulier, les nombres de Hodge $h^{p,q}(F_r(X))$ et $h^{p,q}(G(r,\P^n))$ sont \'egaux si
$p+q<\delta_-$. Rappelons que ces derniers sont nuls pour $p\neq q$,  et qu'ils sont \'egaux si $p=q$
au nombre de partitions de $p$ inscrites dans un rectangle de c\^ot\'es $r+1$ et $n-r$. On retrouve
aussi un r\'esultat de [BV]:

\th 
\label{pic}
Corollaire
\enonce Si de plus $\delta_-\ge 3$, le
groupe de Picard de
  $F_r(X)$ est de rang $1$. \endth

\ex{Remarques} 1) La borne du \theo\ est souvent la meilleure possible: pour une hypersurface 
cubique lisse $X$ dans $\P^5$ et $r=1$, on a pour  
$\delta=4$ et
$\delta_-=2$, et $h^2(F_1(X),\Q )=23$ ([BD], prop. 3); pour une intersection compl\`ete
g\'en\'erique $X$ de deux quadriques dans $\P^6$, on a $h^2(F_1(X),\Q )=8$ ([B2], th.
2.1).\label{rem}

 2) Le \theo\ permet de
retrouver, lorsque  $k=\C$, les points b) et c) du \theo\ \ref{fano}; c'est la m\'ethode suivie
dans [B1].

3)  Il r\'esulte du corollaire \ref{sepunir} et de [K], cor. 1.11, p. 189 et cor.
3.8, p. 202, que pour $n\ge {{\bf d}+r\choose
r+1}+r+1$, les groupes $H^0(F,\Omega^m_F)$ et $H^0(F,(\Omega^1_F)^{\otimes
m})$ sont nuls pour tout $m>0$. Lorsque $k$ est de caract\'eristique nulle, l'annulation de ces
groupes peut se d\'eduire de la remarque \ref{rempi}.2) et du \theo\ 2.13 de [K], p. 254 (\cf\ 
[C2] et [KMM2]), sous l'hypoth\`ese plus faible
$n\ge {{\bf d}+r\choose r+1}$.

4) Lorsque $n\ge {{\bf d}+r\choose
r+1}$,  $F_r(X)$  est une vari\'et\'e de Fano (remarque \ref{rempi}.2). Lorsque $k$ est
de caract\'eristique nulle et que $F_r(X)$ est lisse, le \theo\ d'annulation de Kodaira entra\^ine
que son groupe de Picard est sans torsion (\cf\ [K], (1.4.13), p. 242). Vue l'hypoth\`ese sur $n$, et
sauf dans le cas  
${\bf d}=(2,2)$ et $n=2r+4$, on a  
$\delta_-\ge 3$,
d'o\`u
$\Pic (F_r(X))\isom\Z$ par le corollaire (comparer avec la prop. \ref{pij}.b). Lorsque
${\bf d}=(2,2)$ et $n=2g+1$, la vari\'et\'e $F_{g-2}(X)$ est isomorphe \`a l'espace de  modules des
fibr\'es stables de rang $2$ et de d\'eterminant fix\'e de degr\'e impair sur une courbe
hyperelliptique
$C$ de genre
$g$ ([DR], th. 1). L'isomorphisme
$\Pic (F_{g-2}(X))\isom\Z[\cO(1)]$ (on a  
$\delta_-=3$) est d\'emontr\'e dans [DR], (5.10) (II), p. 177 (\cf\ aussi  [R]). En revanche,
$F_{g-1}(X)$ est isomorphe \`a la jacobienne de $C$ ([DR], th. 2).

\medskip
{\it D\'emonstration du \theo\ \ref{leff}}. Sous les hypoth\`eses du \theo , $F_r(X)$ est le lieu des
z\'eros d'une section du fibr\'e $\Sym ^{\bf d}S^*$, et sa codimension dans la grassmannienne 
est le rang  de ce fibr\'e. Il existe donc une suite exacte (complexe de Koszul):
 $$0\ra\we^{\max}(\Sym^{\bf d}S)\ra\cdots\ra \we^2(\Sym^{\bf d}S) \ra \Sym^{\bf d}S \ra {\cal
I}_{F_r(X)}\ra 0\ .
\leqno{\formule}\hbox{\label{kos}}$$
\ind Notre outil
essentiel sera le \theo\ d'annulation suivant:

\th 
\label{annul}
Proposition
\enonce
Soient $a,b,i,j_1,\ldots ,j_s$ des entiers tels que
$b<a+d_1j_1+\cdots +d_sj_s$ et  $b+i<\delta_-$. Alors
  $$H^{j_1+\cdots +j_s+i}(G(r,\P^n_{\bf C}),\we^{j_1}(\Sym^{d_1}S)\otc \we^{j_s}(\Sym^{d_s}S)
\ot S^{\ot a}\ot S^{*\ot b})=0\ .$$
\endth

\ind Soit $\Sym_{\l}S$ une composante de
  $\we^{j_1}(\Sym^{d_1}S)\otc \we^{j_s}(\Sym^{d_s}S)\ot S^{\ot a}\ot S^{*\ot b}$,
  o\`u\break $\l=(\l_0,\ldots,\l_r)$ est une suite d\'ecroissante d'entiers
  relatifs.  D'apr\`es le th\'eor\`eme de Bott ([De], [Ma1]),
  $H^{j+i}(G,\Sym_{\l}S)$ ne peut \^etre non nul que s'il existe un entier $h$,
  avec $0\leq h\leq r+1$, tel que $j+i=h(n-r)$ et $\l_h\leq h$, ce qui
  implique en particulier que la somme des composantes de $\l$
  d'indice sup\'erieur ou \'egal \`a $h$ v\'erifie $|\l|_{\geq h}\leq
  h(r+1-h)$.

 \ind Comme $|\l|=|\l|_{\geq 0}=d_1j_1+\cdots +d_sj_s+a-b>0$, le cas $h=0$
  est exclu. De plus,
  $$|\l|_{\geq h}\geq j_1+\cdots +j_s-{{\bf d}+h-1 \choose h-1 }-b\ .$$
 \ind En effet, supposons tout d'abord $a=b=0$. Admettons provisoirement le cas $s=1$; le cas
o\`u $s$ est quelconque s'ensuit, puisque si $\Sym_{\l}S$ est un
  facteur direct de $\Sym_{\l_1}S\otc \Sym_{\l_s}S$, la r\`egle de Littlewood
  et Richardson implique  $|\l|_{\ge h}\ge$\break$ |\l_1|_{\ge h}+\cdots
  +|\l_s|_{\ge h}$. Enfin, tensoriser par $S^{ \ot a}$ ne peut qu'augmenter $|\l|_{\ge h}$,
tandis que tensoriser par $S^{*\ot b}$
  fait diminuer $|\l|_{\geq h}$ au plus de $b$. 

 \ind Pour conclure \`a une contradiction, il suffit donc de v\'erifier que
  pour $1\leq h\leq r+1$,
  $$h(n-2r+h-1)-{{\bf d}+h-1 \choose h-1 }>b+i\ .$$
   \ind On retrouve au membre de gauche la fonction  $\phi$ de (\ref{concave}); comme $\delta_-$
est positif, le lemme r\'esulte de l'hypoth\`ese $\delta_->b+i$ comme en (\ref{concave}).
 Il reste \`a traiter le cas $a=b=0$ et $s=1$, qui r\'esulte du lemme suivant. \cqfd

\th
Lemme 
\enonce
 Soient $V$ un espace vectoriel complexe,
$m$ et $d$ des entiers, et $e$ la dimension de $\Sym^dV^m$.
Pour toute composante irr\'eductible $\Sym_{\lambda}V$ de 
$\we^j(\Sym^dV)$, on a $|\lambda |_{>m}\ge j-e$.
\endth

\ind Soit $X$ la grassmannienne des sous-espaces de codimension 
$m$ de $V$, soit $Y$ celle des sous-espaces de codimension $e$ de 
$\Sym^dV$. On notera $S_X$ et $Q_X$ les fibr\'es tautologique et 
quotient sur $X$, de m\^eme que $S_Y$ et $Q_Y$ sur $Y$. 
On plonge $X$ dans $Y$ en associant au noyau du quotient 
$V\rightarrow Q$ celui du quotient induit $\Sym^dV\rightarrow \Sym^dQ$.

\ind D'apr\`es le th\'eor\`eme de Borel-Weil, $\we^j(\Sym^dV)$
est l'espace des sections globales du fibr\'e 
$E=\det Q_Y\ot\we^{j-e}S_Y$. Notons $(\Gamma_l)_{l\ge 0}$ la filtration 
de cet espace de sections selon leur ordre d'annulation $l$ sur $X$.
On dispose d'applications injectives
$$\Gamma_l/\Gamma_{l+1}\hookrightarrow H^0(X,E\vert_X\ot \Sym^lN^*)\ ,$$
o\`u $N$ est le fibr\'e normal de $X$ dans $Y$. 

\ind Le membre de droite 
ne se d\'eduit pas directement du th\'eor\`eme de Borel-Weil. 
Cependant, tout fibr\'e homog\`ene ${\cal F}$ sur $X$ admet une filtration homog\`ene
dont les quotients successifs dont irr\'eductibles, c'est-\`a-dire 
produits de puissances de Schur de $Q_X$ et $S_X$. La somme $\gr {\cal F}$
de ces quotients ne d\'epend pas de la filtration choisie, et le lemme 
de Schur implique l'existence d'une injection $H^0(X,{\cal F})\hookrightarrow
H^0(X,\gr {\cal F})$ : le th\'eor\`eme de Borel-Weil explicite ce dernier
espace de sections. 
Par exemple, $Q_Y\vert_X=\Sym^dQ_X$ est irr\'eductible, et 
$$\gr S_Y\vert_X=\bigoplus_{i=1}^d\Sym^{d-i}Q_X\ot \Sym^iS_X$$
a tous ses termes de degr\'e sup\'erieur ou \'egal \`a $1$ en $S_X$.
Cela implique que $\rm E_{|X}$ est somme de fibr\'es de la 
forme  $\Sym_{\alpha}Q_X\ot \Sym_{\beta}S_X$, avec $|\beta|\ge j-e$.
L'espace des sections globales d'un tel fibr\'e est une puissance de 
Schur $\Sym_{\lambda}V$, o\`u $\lambda=(\alpha,\beta)$ est la partition 
(si c'en est une) obtenue en concat\'enant $\alpha$ et $\beta$. 
En particulier, 
$|\lambda|_{>m}=|\beta|\ge j-e$, ce qui d\'emontre le lemme pour
les composantes de $\we^j(\Sym^dV)$ qui proviennent de $\Gamma_0/\Gamma_1$.

\ind Pour \'etendre ce r\'esultat \`a celles qui proviennent de $\Gamma_l/\Gamma_{l+1}$
pour tout $l>0$, il suffit de s'assurer que toute composante 
irr\'eductible de ${\rm gr} N^*$ est de degr\'e positif ou nul en $S_X$.
Mais c'est une cons\'equence imm\'ediate du fait que $N^*$ est un sous-fibr\'e
homog\`ene de $\Omega^1_Y\vert_X=Q^*_Y\vert_X\ot S_Y\vert_X$, puisque $Q_Y\vert_X$
est de degr\'e z\'ero, et chaque composante de $S_Y\vert_X$ de degr\'e positif
en $S_X$.\cqfd

\ind Revenons \`a la d\'emonstration du \theo\ \ref{leff}; posons $G=G(r,{\bf P}V)$ et $F=F_r(X)$.
Il suffit de le v\'erifier pour la cohomologie complexe, 
  donc, via la d\'ecomposition de Hodge, de d\'emontrer que les morphismes
  $H^q(G,\Omega_G^p)\ra H^q(F,\Omega_F^p)$ sont bijectifs
  pour $p+q<\delta_-$, et injectifs pour
$p+q=\delta_-$. On va montrer que les
  morphismes $H^q(G,\Omega_G^p)\ra H^q(F,\Omega_G^p\vert_F)$ et
  $H^q(F,\Omega_G^p\vert_F)\ra H^q(F,\Omega_F^p)$ ont les m\^emes
  propri\'et\'es.
  
\ind  Pour les premiers, il s'agit de v\'erifier que
  $H^q(G,{\cal I}_F\ot\Omega^p_G)=0$ pour $p+q\leq \delta_-$, donc, via le
  complexe de Koszul, que
  $$H^{q+j-1}(G,\Omega_G^p\ot\we^j(\Sym^{\bf d}S))=0 \quad \hbox{pour tout}\ \  j>0\ .$$
 \ind Rappelons que si $Q$ est le fibr\'e quotient sur $G$, on dispose d'un
  isomorphisme $\Omega^1_G\isom Q^*\ot S$, d'o\`u la suite exacte
  $0\ra\Omega^1_G\ra V^*\ot S \ra S^*\ot S \ra 0$. Sa puissance
  ext\'erieure $p$\tx i\`eme montre que l'annulation pr\'ec\'edente est
  cons\'equence de
  $$H^{q+j-i-1}(G,\we^j(\Sym^{\bf d}S)\ot \we^{p-i}(V^*\ot S)\ot \Sym^i(S^*\ot
  S))=0 \qquad\hbox{pour tout}\ \  j>0\ ,\ i\geq 0\ ,$$
  ce qu'assure la proposition \ref{annul} d\`es que
  $q\leq\delta_-$.
  
\ind   Pour les seconds, la suite exacte normale montre qu'il suffit de
  s'assurer que
  $$H^{q+i}(F,\Omega_G^{p-i-1}\vert_F\ot \Sym^i(\Sym^{\bf d}S))=0 \qquad\hbox{pour tout}\ \  i>0\ ,$$
  donc, \`a cause encore une fois du complexe de Koszul, que
  $$H^{q+i+j}(G,\Omega_{G}^{p-i-1}\ot \Sym^i(\Sym^{\bf d}S)\ot\we^j(\Sym^dS))=0
  \qquad\hbox{pour tout}\ \  i>0\ ,   j\geq 0\ .$$
 \ind En raisonnant comme on vient de le faire, on
  constate que cette annulation a lieu d\`es que $i+q<\delta_-$, ce qui
  conclut cette d\'emonstration puisque $i<p$. \cqfd

\section{Normalit\'e projective, \'equations, degr\'e des sch\'emas de Fano}
 
\th
Th\'eor\`eme
\enonce
 Soit $X$ un sous-sch\'ema de $\P^n_{\bf C}$
d\'efini par des \'equations de degr\'e ${\bf d}$, tel que
 $F_r(X)$  soit de dimension $\delta$. Supposons 
 $n\ge r+{ {\bf d}+r\choose r}$. Alors
$F_r(X)$
  est projectivement normale, autrement dit les morphismes de
  restriction $$\rho_l:H^0(G(r,\P^n),\cO(l))\lra H^0(F_r(X),\cO(l))$$
  sont surjectifs pour tout $l\geq 0$. Par ailleurs, $\rho_l$ est injectif 
pour 
  $l< d_-=\min \{ d_1,\ldots ,d_s\}$. 
\endth

\ind Posons $G=G(r,\P^n)$; d'apr\`es le th\'eor\`eme de Bott,
  $$H^j(G,\we^j(\Sym^{{\bf d}}S)(l))=0\qquad\hbox{pour tout}\ \   j>0\ \ \hbox{et tout}\ \  l\geq 0\
.$$
\ind  En effet, si l'on raisonne comme dans la d\'emonstration de la proposition \ref{annul},
  cet espace ne peut \^etre non nul que si $j$ est multiple de $n-r$; vue l'hypoth\`ese 
$n-r\ge \codim F_r(X)$, la seule possibilit\'e
  est $j=n-r=\codim F_r(X)$, auquel cas $\we^j(\Sym^{\bf d}S)(l)$ est une
  puissance de $\cO(1)$, et n'a donc pas non plus de cohomologie en
  degr\'e $n-r$. La normalit\'e projective s'ensuit, via le complexe de
  Koszul (\ref{kos}) tordu par $\cO(l)$.  

 \ind En fait, les arguments pr\'ec\'edents impliquent plus pr\'ecis\'ement que la
  suite spectrale associ\'ee \`a ce complexe de Koszul tordu d\'eg\'en\`ere en
  $E_2$, ce dont on d\'eduit que le complexe des sections globales
  $$\cdots\lra H^0(G,\we^2(\Sym^{\bf d}S)(l))\lra H^0(G,\Sym^{\bf d}S(l))\lra
  H^0(G,\cI_{F_r(X)}(l))\lra 0$$
  est exact. Mais pour $l<d_-$, on a
  $H^0(G,\Sym^{\bf d}S(l))=0$ d'apr\`es le th\'eor\`eme de Bott, d'o\`u l'inexistence
  d'\'equations de $F_r(X)$ de degr\'e $l$.  \cqfd 

\ex{Remarques} 1) Ce dernier complexe implique au passage que 
$H^0(G,\cI_{F_r(X)}(d_-))$ n'est pas nul, et l'on peut calculer 
explicitement sa dimension.

2) Les sch\'emas de Fano ne sont en g\'en\'eral pas projectivement normaux; si l'on revient au cas
${\bf d}=(2,2)$ et $n=2g+1$ (\cf\ rem. \ref{rem}.4), Laszlo a montr\'e dans [L] (par des
m\'ethodes similaires) que $H^0(F_{g-2},\cO(1))$ est de dimension $2^{g-1}(2^g-1)$. En
particulier, $\rho_1$ n'est pas surjectif.

3) Le
\theo\ d'annulation de Kodaira entra\^ine que les groupes
$H^i(F_r(X),\cO(l))$ sont nuls  pour $i>  0$ et $l\ge -n+{{\bf d}+r\choose
r+1}$. Si l'on raisonne comme dans la preuve de la
proposition
\ref{annul}, 
    on montre facilement la m\^eme annulation pour $i< 
    \min (\delta ,n-(l+2)r-s)$. A l'ext\'erieur du domaine d\'efini par ces in\'egalit\'es,
    il peut ne pas y avoir annulation: pour une sextique $X$ dans $\P^6$,  on peut montrer que $
H^2(F_1(X),\cO (6))$ est de dimension $\ge 10024$ (alors que $F_1(X)$ est de dimension $3$).

\bigskip
\ind Introduisons des polyn\^omes de $r+1$ variables, $e(x)=x_0+\cdots +x_r$,
et
$$Q_{r,d}(x)=\prod_{a_0+\cdots +a_r=d}(a_0x_0+\cdots +a_rx_r)\ ,$$ 
puis $Q_{r,{\bf d}}(x)=Q_{r,d_1}(x)\cdots Q_{r,d_s}(x)$.

\th\label{calcul}
Th\'eor\`eme
\enonce
 Soit $X$ un sous-sch\'ema de $\P^n_k$
d\'efini par des \'equations de degr\'e ${\bf d}$, tel que
 $F_r(X)$  soit de dimension $\delta$.
 Le degr\'e de $F_r(X)$ pour le plongement de Pl\"ucker de $G(r,\P^n)$ 
  est \'egal au coefficient du mon\^ome $x_0^nx_1^{n-1}\cdots x_r^{n-r}$
  dans le produit du polyn\^ome $Q_{r,{\bf d}}\times e^{\delta}$
  et du Vandermonde.
\endth

\ind Ce degr\'e est
  $$\deg (F)=\int_{G(r,{\bf P}^n)} c_{\max}(\Sym^{\bf d}S^*)c_1(\cO(1))^{\delta}\ .$$
 \ind Rappelons
  que l'anneau de Chow de $G(r,\P^n)$ est un quotient d'un
  anneau des polyn\^omes sym\'etriques de $r+1$ variables $x_0,\ldots
  ,x_r$, $e(x)$ relevant $c_1(\cO(1))$, et $Q(x)$ relevant
  $c_{\max}(\Sym^dS^*)$ ([Fu], 14.7).  De plus, int\'egrer sur $G$ revient, au
  niveau des polyn\^omes, \`a d\'ecomposer sur les polyn\^omes de Schur
([M]), et ne retenir que le coefficient de celui qui est
  associ\'e \`a la partition rectangle ayant $r+1$ parts \'egales \`a $n-r$, \`a
  savoir $(x_0\ldots x_r)^{n-r}$.

\ind  Il suffit donc de montrer que si $P$ est un polyn\^ome sym\'etrique, que
  l'on d\'ecompose sur les polyn\^omes de Schur, le coefficient du
  pr\'ec\'edent est \'egal \`a celui du mon\^ome $x_0^nx_1^{n-1}\cdots
  x_r^{n-r}$ dans le produit de $P$ et du Vandermonde. Mais par
  lin\'earit\'e, il suffit de le v\'erifier lorsque $P$ est lui-m\^eme un
  polyn\^ome de Schur, auquel cas c'est une cons\'equence imm\'ediate de
  l'expression de Jacobi de ces polyn\^omes ([FH], (A.23), p. 459). \cqfd

\ind Donnons quelques exemples num\'eriques, d'abord pour le cas des droites
d'une hypersurface, qui est d\^u \`a Van der Waerden ([vW]), puis pour
$r\geq 2$, toujours dans le cas d'une hypersurface.

{\eightpoint
$$\vbox{\offinterlineskip \halign{\tv\hskip.5mm \tv#&\cc{$#$}&  \tv#&
\cc{$#$}&  \tv#&\cc{$#$}&  \tv#&\cc{$#$}& \tv\hskip.5mm 
\tv#&\cc{$#$}& \tv#&\cc{$#$}& \tv#&\cc{$#$}& 
\tv#&\cc{$#$}&  \tv\hskip.5mm \tv#\cr
 \noalign{\hrule}\tvi
 &d&& n  &&  \dim F && \deg F  &&d&& n  &&  \dim F && \deg F &\cr
 \noalign{\hrule}\tvi
 &3&& 3 &&  0&&     27 &&  5&& 5 &&  2&&    6\;125&\cr
 \noalign{\hrule}\tvi
 &3&& 4 &&  2&&     45 &&  5&& 6 &&  4&&   16\;100&\cr
  \noalign{\hrule}\tvi
&3&& 5 &&  4&&    108 &&  5&& 7 &&  6&&   46\;625&\cr
 \noalign{\hrule}\tvi
 &4&& 4 &&  1&&    320 &&  6&& 5 &&  1&&  60\;480&\cr
 \noalign{\hrule}\tvi
 &4&& 5 &&  3&&    736 &&  6&& 6 &&  3&& 154\;224&\cr
 \noalign{\hrule}\tvi
 &4&& 6 &&  5&& 1\;984 &&  7&& 5 &&  0&& 698\;005&\cr
 \noalign{\hrule}\tvi
 &4&& 7 &&  7&& 5\;824 &&  7&& 6 &&  2&& 1\;707\;797&\cr
 \noalign{\hrule}\tvi
 &5&& 4 &&  0&& 2\;875 &&  9&& 6 &&  0&& 305\;093\;061&\cr
\noalign{\hrule} }}$$}
\centerline{{\eightpoint 1. Degr\'es de sch\'emas de Fano de droites ($r=1$).}}

\bigskip
{\eightpoint
$$\vbox{\offinterlineskip \halign{\tv\hskip.5mm \tv#&\cc{$#$}&  \tv#&
\cc{$#$}&  \tv#&
\cc{$#$}&  \tv#&\cc{$#$}&  \tv#&\cc{$#$}& \tv\hskip.5mm 
\tv#&\cc{$#$}& \tv#&\cc{$#$}&  \tv#&
\cc{$#$}& \tv#&\cc{$#$}& 
\tv#&\cc{$#$}&  \tv\hskip.5mm \tv#\cr
 \noalign{\hrule}\tvi
 &r &&d&& n  &&  \dim F && \deg F  &&r&&d&& n  &&  \dim F && \deg F &\cr
 \noalign{\hrule}\tvi
 &2&& 3 &&  6&&   2 &&       2\;835&& 2 &&  5 && 9&& 0 && 2\;103\;798\;896\;875&\cr
 \noalign{\hrule}\tvi
 &2&& 3 &&  7&&   5 &&      24\;219&& 3 &&  3 && 8&& 0 &&      321\;489&\cr
 \noalign{\hrule}\tvi
 &2&& 3 &&  8&&   8 &&     274\;590&& 3 &&  3 && 9&& 4 &&   10\;344\;510&\cr
 \noalign{\hrule}\tvi
 &2&& 4 &&  7&&   0 &&  3\;297\;280&& 4 &&  3 && 11&& 0 &&  1\;812\;768\;336&\cr
\noalign{\hrule} }}$$}
\centerline{{\eightpoint 2. Degr\'es de sch\'emas de Fano pour $r=2,3,4$.}}

\medskip

\ind La m\^eme m\'ethode permet en fait de d\'eterminer la d\'ecomposition
$$[F_r(X)]=\sum_{|\l|=\codim F_r(X)}f_{\l}\s_{\l}$$
de la classe fondamentale de
$F_r(X)$ sur les classes des cycles de Schubert de la grassmannienne, o\`u
l'on note $\s_{\l}$ la classe du cycle de codimension $|\l|$ associ\'e \`a
la partition $\l=(\l_0,\ldots ,\l_r)$.

\th
Proposition
\enonce
Si l'on \'ecrit $Q_{r,{\bf d}}(x)=\sum_{\a}q_{\a}x^{\a}$, 
et si $\k$  d\'esigne la suite $(r,\ldots ,1,0)$, alors
  $$f_{\l}=\sum_{\s\in\cS_{r+1}}\e(\s)q_{\s(\l+\k)-\k}\ .$$
\endth

\ind Notons que si l'on adopte pour les cycles de Schubert la m\^eme
convention que pour les polyn\^omes de Schur, \`a savoir que
pour chaque suite d'entiers $\a$, on pose $\s_{\a}=\e(\t)\s_{\l}$ s'il existe
une partition $\l$ et une permutation $\t\in\cS_{r+1}$ telles que
$\a+\k=\t(\l+\k)$, et $\s_{\a}=0$ sinon, la proposition pr\'ec\'edente
se traduit par la simple \'egalit\'e
$$[F_r(X)]=\sum_{\a}q_{\a}\s_{\a}.$$

\ind Donnons par exemple les classes de quelques vari\'et\'es de Fano en bas degr\'e.

$$\eqalign{
{\rm Si}\ \  {\bf d}=(2) \ ,\qquad [F_r] &= 2^{r+1}\sigma_{r+1,r,\ldots,1}\ ,\cr
{\rm Si}\ \  {\bf d}=(3) \ ,\qquad [F_1] &= 9(2\s_{3,1}+3\s_{2,2})\ ,\cr
             [F_2] &= 27(8\s_{6,3,1}+12\s_{6,2,2}+20\s_{5,4,1}
                     +50\s_{5,3,2}+42\s_{4,4,2}+35\s_{4,3,3})\ .\cr
{\rm Si}\ \   {\bf d}=(4)\  ,\qquad  [F_1] &= 32(3\s_{4,1}+10\s_{3,2})\ ,\cr
             [F_2] &= 512(54\s_{10,4,1}+180\s_{10,3,2}+369\s_{9,5,1}+
                     1599\s_{9,4,2}+1230\s_{9,3,3}\cr
&\qquad +819\s_{8,6,1}+5022\s_{8,5,2}+8459\s_{8,4,3}+504\s_{7,7,1}+
                     6039\s_{7,6,2}\cr
&\ \ +18889\s_{7,5,3}+13354\s_{7,4,4}+
                     11660\s_{6,6,3}+23560\s_{6,5,4}+6440\s_{5,5,5})\ .\cr
{\rm Si}\ \   {\bf d}=(5)\  ,\qquad [F_1] &= 25(24\s_{5,1}+130\s_{4,2}+91\s_{3,3})\ .\cr
{\rm Si}\ \   {\bf d}=(2,2)\ ,\quad [F_1] &= 16(\s_{4,2}+\s_{3,3})\ ,\cr
             [F_2] &= 64(\s_{6,4,2}+\s_{6,3,3} +\s_{5,5,2}+2\s_{5,4,3}
                    +\s_{4,4,4})\ .\cr}$$

 \section {Espaces \lin s sur les sch\'emas de Fano} 
 
\ind Le but de ce paragraphe est de montrer que les sch\'emas de Fano sont
s\'eparablement unir\'egl\'es en droites pour
$n$ assez grand (corollaire \ref{sepunir}). Pour cela, nous commen\c cons par g\'en\'eraliser
les r\'esultats du \S 2 aux sous-sch\'emas de $F_r(X)$ form\'es  des $r$\tx plans contenant un
sous-espace
\lin\ fixe de dimension $r_0<r$. Pour de tels entiers, on pose  
$$ \delta(n,{\bf d},r,r_0)=(r-r_0)(n-r)+{{\bf d}+r_0\choose r_0}- {{\bf d}+r\choose r}
$$
et $$\delta_-(n,{\bf d},r,r_0)=\min\{ \delta (n,{\bf d},r,r_0), n-2r+r_0+1-{{\bf d}+r_0\choose
r_0+1}\}\ ,$$de sorte que $\delta(n,{\bf d},r)=\delta_-(n,{\bf d},r,-1)$ et $\delta_-(n,{\bf
d},r)=\delta_-(n,{\bf d},r,-1)$. De nouveau, il est utile de noter que  lorsque ${\bf d}\ne
(2)$, l'entier
  $\delta(n,{\bf d},r,r_0)$ est positif (resp. strictement positif)  si et seulement si
$\delta_-(n,{\bf d},r,r_0)$ l'est; cela
r\'esulte de la convexit\'e de la fonction $\psi:r\mapsto {{\bf d}+r\choose
r}-r^2$, qui entra\^ine l'in\'egalit\'e $\psi(r)-\psi(r_0)\ge
(r-r_0)(\psi(r_0+1)-\psi(r_0))$ (puisque $r>r_0$). Le \theo\ suivant g\'en\'eralise le \theo\
\ref{fano}.
  
\th
\label{fano2}
Th\'eor\`eme
\enonce
Soit $X$ un sous-sch\'ema de
$\P^n_k$ d\'efini par des \'equations de degr\'e ${\bf d}$, soit $\Lambda_0$ un  $r_0$\tx plan contenu
dans $X$, et  soit $F_r(X,\Lambda_0)$, avec $r> r_0$, le sch\'ema de Hilbert des $r$\tx plans
contenus dans
$X$
 et contenant $\Lambda_0$.

\indp {\rm a)} Lorsque $\delta_-(n,{\bf d},r,r_0)< 0$, le sch\'ema $F_r(X,\Lambda_0)$ est vide pour
$X$ g\'en\'erale et $\Lambda_0$ g\'en\'eral contenu
dans $X$.

\indp {\rm b)} Lorsque $\delta_-(n,{\bf d},r,r_0)\ge 0$, le sch\'ema $F_r(X,\Lambda_0)$ est non vide;
il est lisse de dimension
$\delta(n,{\bf d},r,r_0)$ pour $X$ g\'en\'erale et $\Lambda_0$ g\'en\'eral contenu
dans $X$.

\indp {\rm c)} Lorsque $\delta_-(n,{\bf d},r,r_0)>0$, le sch\'ema $F_r(X,\Lambda_0)$ est connexe.
\endth

\ind En gardant les notations de la d\'emonstration du \theo\ \ref{fano},  on
consid\`ere\break $G_0=\{ [\Lambda]\in G(r,{\bf P}^n)\mid \Lambda\supset
\Lambda_0\}$. La dimension de
$I_0=q^{-1}(G_0)$ est \'egale \`a 
$$\dim \P
\Sym^{\bf d}V^*-{{\bf d}+r\choose r}+(r-r_0)(n-r)\ .$$\ind Le c\^one $S_0$ dans
$\Sym^{\bf d}V^*$
 correspondant aux sous-sch\'emas contenant
$\Lambda_0$ est de codimension ${{\bf d}+r_0\choose r_0}$, de sorte
que 
$\dim I_0 =\dim  \P S_0 +\delta$. Supposons $\d_-<0$; si ${\bf d}=(2)$, cela signifie $2r\ge n$, et
on a d\'ej\`a vu qu'une quadrique lisse dans $\P^n$ ne contenait pas de $r$\tx plan; si ${\bf d}\ne
(2)$, on a $\d<0$, et le morphisme $p_0:I_0\ra \P S_0$ induit par $p$ n'est pas surjectif. 

\ind Cela montre a); on
suppose maintenant $\d_-\ge 0$.
Fixons un $r$\tx plan $\Lambda$ contenant $\Lambda_0$, et choisissons des
coordonn\'ees de fa\c con que $\Lambda_0$ soit d\'efini par les \'equations $x_{r_0+1}=\cdots=x_n=0$,
et $\Lambda$ par  $x_{r+1}=\cdots=x_n=0$; pour tout entier positif $m$, on note
$\Gamma_0(m)$ le noyau du morphisme $\Gamma_\Lambda (m)\ra \Gamma_{\Lambda_0} (m)$. 
 
\ind La d\'emarche est enti\`erement analogue \`a celle de la d\'emonstration de \ref{fano}.
Soit ${\bf f}$ un \'el\'ement de $S_0$; pour que $p_0$ soit lisse en un point
$(X_{\bf f},\Lambda)$ de
$I_0$, il faut et il suffit que le morphisme
$\alpha_0:\Gamma_0(1)^{n-r}\ra \Gamma_0({\bf d})$ induit par le morphisme
$\alpha$ du lemme \ref{lisse} soit surjectif.

\ind Soit $Z_0$ le lieu des points de $I_0$ o\`u $p_0$ n'est pas lisse; on montre comme en
\ref{propre}--\ref{ferme}, par r\'ecurrence sur $r-r_0$, que $p_0(Z_0)$ est distinct de
$\P S_0$. Soit
$\mu_0: \Gamma_0(1)\times \Gamma_\Lambda({\bf d}-1)\ra \Gamma_0({\bf d})$
le morphisme induit par la multiplication $\mu$. On montre de la m\^eme fa\c con que
si
$h$ est un entier compris entre $1$ et $r-r_0$, l'ensemble
 des formes lin\'eaires $\ell_0$ sur
$\Gamma_0({\bf d})$
 telles que
$\codim  \mu_0^{-1}(\ell_0)  =h$ est de dimension 
$$ \le
h(r-r_0-h)+{{\bf d}+r_0+h\choose r_0+h}-{{\bf d}+r_0\choose r_0}\ .$$      
\ind On en d\'eduit que la codimension de $Z_0$ dans $I_0$ est 
$$\eqalign{&\ge \min_{1\le h\le
r-r_0}[h(n-r)-h(r-r_0-h)-{{\bf d}+r_0+h\choose r_0+h}+{{\bf d}+r_0\choose r_0}]+1\cr
&=\min\{ n-2r+r_0+1-{{\bf d}+r_0\choose r_0+1}, \delta\}+1=\d_-+1\ ,\cr}$$
puisque la fonction entre crochets est une fonction concave de $h$ lorsque ${\bf d}\ne (2)$, et
croissante lorsque ${\bf d}=(2)$ puisque $\d_-$ est positif (\cf\ (\ref{concave})). La fin de la
d\'emonstration est la m\^eme que celle du \theo\ \ref{fano}. \cqfd

\ind Soient
 $X$ un sous-sch\'ema de
$\P^n_k$ d\'efini par des \'equations de degr\'e ${\bf d}$, et $\Lambda$ un $(r+1)$\tx plan
contenu dans
$X$. Les $r$-plans contenus dans $\Lambda$ d\'efinissent une inclusion
de $\Lambda^*$ dans
$F_r(X)$, dont l'image par le plongement de Pl\"ucker est  un $(r+1)$\tx plan.

\th
\label{uni}
Corollaire
\enonce
Soit $X$ un sous-sch\'ema de
$\P^n_k$ d\'efini par des \'equations de degr\'e ${\bf d}$.

\indp {\rm a)} Si ${\bf d}\ne 2$ et $n\ge {1\over r}{{\bf d}+r\choose r}+r-{s\over r}$, ou si
${\bf d}= 2$ et $n\ge 2r+1$, la
vari\'et\'e $X$ est recouverte par des $r$\tx plans.

\indp {\rm b)}  Si $n\ge {{\bf d}+r\choose r+1}+r+1$, la sous-vari\'et\'e 
$F_r(X)$ de $G(r,\P^n)$ est unir\'egl\'ee en droites.
\endth

\ind Le point a) r\'esulte du \theo\ avec $r_0=0$. Soit $\Lambda_0$ un $r$\tx plan
contenu dans
$X$; sous les hypoth\`eses de b),  le \theo\ \ref{fano2}.b)  entra\^ine qu'il
existe un
$(r+1)$\tx plan
 $\Lambda_1$ contenu dans $X$ et contenant $\Lambda_0$. Le $(r+1)$\tx plan $\Lambda_1^*$,
contenu dans $F_r(X)$, passe par $[\Lambda_0]$. En particulier, il passe une droite par tout
point de
$F_r(X)$.\cqfd

\th
Th\'eor\`eme
\enonce
Soit $X$ un sous-sch\'ema g\'en\'eral de
$\P^n_k$ d\'efini par des \'equations de degr\'e ${\bf d}$; on suppose  $n\ge {{\bf d}+r\choose
r+1}+r+1$. Soit $\Lambda$ un
$(r+1)$\tx plan g\'en\'eral contenu dans $X$. La restriction \`a une droite g\'en\'erale de
$\Lambda^*$ du fibr\'e normal \`a
$\Lambda^*$ dans $F_r(X)$  est isomorphe \`a
$$ {\cal O}^{r(n-r-1)+{{\bf d}+r\choose r+1}-{{\bf d}+r\choose r}}\oplus{\cal
O}(1)^{n-r-1 -{{\bf d}+r\choose r+1}}\ .$$ 
\endth

\ind Soit $N$ le fibr\'e normal \`a
$\Lambda^*$ dans $F_r(X)$; on a la suite exacte
\vskip -5mm
$$\diagram{ 0&\lra &N  &\lra &N_{\Lambda^*/G} &\lra
&\bigl( N_{F_r(X)/G}\bigr) \vert_{\Lambda^*} &\lra &0\cr
&&&&||&&||\cr
&&&&(S^*\vert_{\Lambda^*})^{ n-r-1 }&&\Sym^dS^*\vert_{\Lambda^*}\cr}$$
\vskip -5mm
 dont la restriction \`a une droite $\ell$ contenue dans
$\Lambda^*$ est
$$ 0  \lra  N \vert_\ell  \lra   ( S^*\vert_\ell)^{ n-r-1 } 
\buildrel{u}\over{\lra}  \Sym^d S^*\vert_\ell\lra  0\ . \leqno{\formule}\hbox{\label{suite}}$$
\ind Comme $S^*\vert_\ell$ est isomorphe \`a $\cO ^r\oplus\cO (1)$, cela entra\^ine que
$N_\ell$ est isomorphe
\`a une somme directe 
$\bigoplus_j\cO(a_j)$ avec $a_j\le 1$ pour tout $j$. On v\'erifie que
$H^0(\ell,S^*\vert_\ell)$ s'identifie
\`a
$H^0(\Lambda,\cO (1))$,
\cad\  
\`a l'espace vectoriel not\'e $\Gamma_{\Lambda}(1)$ dans la d\'emonstration du \theo\
\ref{fano} et $H^0(\ell,\Sym^dS^*\vert_\ell)$  \`a $\Gamma_{\Lambda}(d)$. Soient $x_0$ un point
de $\ell$, et $\Lambda_0$ l'hyperplan de $\Lambda$ associ\'e. On a un diagramme commutatif
$$\diagram{\Gamma_0(1)^{n-r-1}
&\phfl{}{}&\Gamma_\Lambda(1)^{n-r-1}&\phfl{}{}&\Gamma_{\Lambda_0}(1)^{n-r}\cr
\pvfl{\displaystyle\alpha_0}{}&&\pvfl{\displaystyle\alpha}{}&&\pvfl{}{}\cr
\Gamma_0(d)
&\phfl{}{}&\Gamma_\Lambda(d)&\phfl{}{}&\Gamma_{\Lambda_0}(d)\cr
}$$  
o\`u les notations sont celles de la d\'emonstration du \theo\ \ref{fano2}. On v\'erifie que
$\alpha$ s'identifie \`a $H^0(u)$, et
$\alpha_0$ \`a
$H^0(u(-x_0)):H^0(\ell, ( S^*\vert_\ell)(-x_0)^{ n-r-1 } )\ra H^0(\ell, \Sym^d
S^*\vert_\ell(-x_0) )$. Comme 
$$\delta_-(n,{\bf d}, r+1,r)= n-r-1+{{\bf d}+r\choose r+1} $$
est positif par hypoth\`ese, la d\'emonstration du \theo\ \ref{fano2} entra\^ine que 
$H^0(u(-x_0))$ est surjectif; il en r\'esulte que $H^1(\ell,N\vert_\ell(-x_0))$ est nul.
Cela entra\^ine que les $a_j$ sont tous positifs. Le rang et le degr\'e de
$N_\ell$
\'etant donn\'es par (\ref{suite}), cela d\'emontre le \theo .\cqfd

\ind Il n'est pas vrai en g\'en\'eral que le fibr\'e normal \`a
$\Lambda^*$ dans $F_r(X)$ soit somme de fibr\'es en droites; cependant,
c'est le cas lorsque
$\delta(n,{\bf d},r+1)$ est nul ([BV], prop. 3). 
  
\th
\label{sepunir}
Corollaire
\enonce
Soit $X$ un sous-sch\'ema g\'en\'eral de
$\P^n_k$ d\'efini par des \'equations de degr\'e ${\bf d}$; on suppose  $n\ge {{\bf d}+r\choose
r+1}+r+1$. La vari\'et\'e $F_r(X)$
est s\'eparablement unir\'egl\'ee en droites. 
\endth

\ind L'hypoth\`ese sur $n$ entra\^ine que $\delta_-(n,{\bf d},r+1)$ est
positif; soient $\Lambda_1$ un $(r+1)$\tx plan g\'en\'eral contenu dans
$X$, et $\ell$ une droite g\'en\'erale contenue dans $\Lambda_1^*$. Le
\theo\ pr\'ec\'edent entra\^ine que le fibr\'e normal \`a $\ell$ dans 
$F_r(X)$ est somme de copies de $\cO_\ell$ et $\cO_\ell(1)$, donc que
$\ell$ est {\it libre} au sens de [K], p. 113 (et m\^eme  {\it minimale}
au sens de \loc , p. 195). Le corollaire r\'esulte alors de
\loc , p. 188.\cqfd

\section{Cycles alg\'ebriques}

\ind On voudrait montrer que pour $n$ assez grand, les groupes de Chow
rationnels de $F_r(X)$ sont les m\^emes que ceux de la grassmannienne
ambiante $G(r,\P^n)$, g\'en\'eralisant ainsi des r\'esultats de [P], [K]
p. 266, et  [ELV], qui traitent le cas $r=0$. On n'obtient
malheureusement  de r\'esultats nouveaux que pour les groupes
$A_1(F_r(X))_{\bf Q}$, en caract\'eristique nulle. Les id\'ees sont celles
de [K].

\th
\label{rat}
Proposition
\enonce
Soit $X$ un sous-sch\'ema de
$\P^n_k$ d\'efini par des \'equations de degr\'e ${\bf d}$; on suppose
$n\ge {{\bf d}+r\choose r+1}$. Le sch\'ema 
$F_r(X)$ est connexe par cha\^ines rationnelles; en
particulier, $A_0(F_r(X))\isom\Z$. 
\endth

\ind Lorsque $X$ est g\'en\'erale, il r\'esulte du  \theo\
\ref{fano} et de la remarque  \ref{rem}.4) que $F_r(X)$ est une vari\'et\'e de Fano  lisse
connexe, donc est connexe par cha\^ines rationnelles ([K], 2.13, p. 254).
Le cas g\'en\'eral s'en d\'eduit comme dans [K], 4.9, p. 271.\cqfd

\ind On suppose maintenant $k=\C$ (pour g\'en\'eraliser les
r\'esultats qui suivent  en toute caract\'eristique, il suffirait de montrer que le groupe
de N\'eron-Severi d'un sch\'ema de Fano g\'en\'eral est de rang $1$).

\th 
\label{chaine}
Proposition
\enonce
Soit $X$ un sous-sch\'ema de
$\P^n_{\bf C}$ d\'efini par des \'equations de degr\'e ${\bf d}$; on suppose $n\ge {{\bf
d}+r\choose r+1}+r+1$. Deux points quelconques de $F_r(X)$ peuvent \^etre joints par une
courbe connexe de degr\'e
$\delta(n,{\bf d},r)$, dont toutes les composantes sont des droites.
\endth

\ind On peut supposer $X$ g\'en\'erale, de sorte que $F_r(X)$ est une
vari\'et\'e de Fano lisse  unir\'egl\'ee en droites (cor. \ref{uni}.b)), de groupe de
N\'eron-Severi de rang $1$ (cor.  \ref{pic}). Le corollaire r\'esulte de [K], 
 p. 252.\cqfd

\ind Soient $X$ un $k$\tx sch\'ema et $m$ un entier positif;  on
 note   
$A_m(X)$ (resp. $B_m(X)$) le groupe des classes d'\'equivalence rationnelle (resp.
alg\'ebrique) de
$m$\tx cycles sur
$X$ (\cf\ [K], p. 122).

\th 
Th\'eor\`eme
\enonce
Soit $X$ un sous-sch\'ema de
$\P^n_{\bf C}$ d\'efini par des \'equations de degr\'e ${\bf d}$.

\indp {\rm a)} Si $n\ge {{\bf d}+r\choose r+1}+r+1$, l'espace vectoriel
$B_1(F_r(X))_{\bf Q}$ est de rang $1$.

\indp {\rm b)} Si $n\ge {{\bf d}+r+1\choose r+2}$, l'espace vectoriel
$A_1(F_r(X))_{\bf Q}$ est de rang $1$.
\endth

\ind En utilisant IV.3.13.3 de [K], p. 206, et en raisonnant comme dans \loc , p. 271, on voit
que le corollaire \ref{chaine} entra\^ine que
$A_1(F_r(X))_{\bf Q}$ est engendr\'e par les classes des droites. Ces droites sont param\'etr\'ees
par un fibr\'e en $G(r-1,\P^{r+1})$ au-dessus de $F_{r+1}(X)$, de sorte qu'il existe un
morphisme surjectif  $A_0(F_{r+1}(X))_{\bf Q}\ra A_1(F_r(X))_{\bf Q}$. Sous l'hypoth\`ese de a),
$F_{r+1}(X)$ est connexe. Sous l'hypoth\`ese de b), il r\'esulte du  cor. \ref{rat}
que $A_0(F_{r+1}(X))_{\bf Q}$ est de dimension $1$.\cqfd 

\ind Lorsque $F_r(X)$ est lisse, la conclusion de la partie a) du \theo\ pr\'ec\'edent
reste  valable sous l'hypoth\`ese plus faible
$n\ge {{\bf d}+r\choose r+1}$; cela r\'esulte du corollaire \ref{pic} et de [BS] (\cf\ aussi [K], th.
3.14, p. 207) .

\ind Lorsque  $X$ contient un $(r+l)$\tx
plan $\Lambda$, le plongement $G(r,\Lambda)\i F_r(X)\i G(r,\P^n)$ induit un
isomorphisme $A_i(G(r,\Lambda))\isom A_i(G(r,\P^n))$ pour $i\le l$ ([Fu], p. 271), de sorte
qu'on a une surjection $A_i(F_r(X))\twoheadrightarrow A_i(G(r,\P^n))$.
 
\th 
Conjecture
\enonce
Soit $X$ un sous-sch\'ema de
$\P^n_k$ d\'efini par des \'equations de degr\'e ${\bf d}$. Si $n\ge {{\bf d}+r+l\choose
r+l+1}$, le morphisme 
$A_l(F_r(X))_{\bf Q}\ra A_l(G(r,\P^n))_{\bf Q}$ induit par l'inclusion est bijectif.
\endth

\ind Lorsque $l=1$ et $k=\C$, c'est le \theo\ pr\'ec\'edent; pour $r=0$ c'est le \theo\ principal
de [ELV].

\section{Unirationalit\'e}

\ind Nous allons maintenant d\'emontrer  l'unirationalit\'e de certains sch\'emas de Fano  en
nous ramenant \`a un r\'esultat de [PS], qui fournit un crit\`ere explicite pour
l'unirationa\-li\-t\'e d'une intersection compl\`ete dans un espace projectif. 

\ind Ce crit\`ere est le
suivant. On d\'efinit tout d'abord, pour toute suite ${\bf d}=(d_1,\ldots,d_s)$ d'entiers strictement
positifs, des entiers $n({\bf d})$ et $r({\bf d})$ de la fa\c con suivante: on pose
$n(1)=r(1)=0$ (dans [PS], on trouve $n(1)=1$, mais $n(1)=0$ suffit); si l'un des
$d_i$ vaut $1$, on note ${\bf d'} $ la suite ${\bf d}$ priv\'ee de
$d_i$, et on pose
$n({\bf d})=n({\bf d'})+1$ et $r({\bf d})=r({\bf d'})$; enfin, si tous les $d_i$
sont $>1$, on pose
 $r({\bf d})=n({\bf d }-1)$ et $n({\bf d})=r({\bf d})+{{\bf d}+r({\bf d})-1\choose r({\bf d})}$.
On a par exemple
$$\nospacedmath\displaylines{r(2,\ldots,2)=s-1\qquad\qquad r(3,\ldots,3)=s^2+s-1\cr 
 r(4,\ldots,4 )= s^2+s-1+{s^2(s+1)(s^2+s+1)\over 2}\
.\cr}$$
\ind Les bornes donn\'ees dans [R] sont un peu meilleures, mais je ne sais pas extraire de cet
article un crit\`ere effectif. 

\th\label{PS} 
Th\'eor\`eme ([Pr], [PS])
\enonce
Soit $F$ une intersection compl\`ete dans $\P^N_k$
d\'efinie par  des
\'equations ${\bf
g}=(g_1,\ldots,g_S)$ de degr\'e ${\bf
D}$ et contenant un $r({\bf D})$\tx plan $\Lambda$. On suppose $N\ge n({\bf D})$,
que
$F$ est irr\'eductible de codimension $S$ et lisse le long de $\Lambda$, et que l'application
$\beta:k^{N+1}\ra \Gamma_{\Lambda}({\bf D}-1)$ d\'efinie par
$$\beta(a_0,\ldots,a_N)=\Bigl( \sum_{i=0}^Na_i \Bigl({\partial
g_1\over\partial x_i}\Bigr) _{\displaystyle{\vert_\Lambda}},\ldots,
\sum_{i=0}^Na_i \Bigl({\partial
g_S\over\partial x_i}\Bigr) _{\displaystyle{\vert_\Lambda}}\Bigr)$$
est surjective. Alors $F$ est  unirationnelle.
\endth

\ind On remarquera que la surjectivit\'e de $\beta$ entra\^ine celle de l'application $\alpha$
d\'efinie en \ref{lisse}, donc la lissit\'e de $F_{r({\bf D})}(F)$ en $\Lambda$.

\th
\label{unirat}
Th\'eor\`eme
\enonce
Il existe une constante explicite $n({\bf d},r)$ telle que, pour $n\ge n({\bf d},r)$, le sch\'ema
de Fano des
$r$\tx plans contenus dans un sous-sch\'ema g\'en\'erique $X$
de $\P^n_k$ d\'efini par des \'equations de degr\'e ${\bf d}$, 
 est unirationnel.
\endth

\ex{Remarques} 1) La borne $n({\bf d},r)$ que l'on obtient est tr\`es grande. Elle est d\'efinie de
la fa\c con suivante: soit ${\bf D}$ la suite d'entiers
 o\`u chaque
$d_i$ est r\'ep\'et\'e ${d_i +r\choose r}$ fois; on pose $r_1=(r({\bf D})+1)(r+1)-1$ et
$$n({\bf d},r)= r_1+{{\bf d}+r_1-1\choose r_1}\ .$$

\ind Pour le cas le plus simple $r=1$ et ${\bf d}=(3)$, \cad\ pour le sch\'ema des droites contenues
dans une hypersurface cubique, on a ${\bf D}=(3,3,3,3)$, $r(3,3,3,3)=19$ et $n((3),1)=859$. Dans ce
cas pr\'ecis, il est facile d'am\'eliorer la borne de [PS] en $r(3,3,3,3)=13$ (il suffit de
remarquer qu'une intersection de $4$ quadriques est rationnelle d\`es qu'elle contient un $3$\tx
plan dans son lieu lisse, en proc\'edant par exemple comme dans [CTSSD]); on obtient alors
$n((3),1)=433$.

\ind On obtient aussi $n((2,\ldots,2),r)=s(s+1){r+2\choose 2}(r+1)-1$. Rappelons
que pour ${\bf d}=(2,2)$ et $n=2g+1$, la vari\'et\'e $F_r(X)$ est une \va\ pour $r=g-1$ (\cf\ rem.
\ref{rem}.4), et qu'elle est  {\it rationnelle} pour $r=g-2$ ([N]), donc unirationnelle pour $r\le
g-2$.

2) L'adjectif \og g\'en\'erique\fg\ de l'\'enonc\'e du \theo\ peut \^etre pr\'ecis\'e: si
$n\ge n({\bf d},r)$, le sch\'ema $F_r(X)$ est unirationnel s'il est irr\'eductible de dimension
$\d (n,{\bf d},r)$, si $X$ contient un $r_1$\tx plan $\Lambda_1$ pour lequel l'application $\beta$ du
\theo\ \ref{PS} est surjective, et si $F_r(X)$ est lisse le long de $G(r,\Lambda_1)$.

{\it D\'emonstration du \theo }. Soit $V$ l'espace vectoriel $k^{n+1}$. On note
$(x^{(0)},\ldots,x^{(r)})$, avec
$x^{(j)}=(x^{(j)}_0,\ldots,x^{(j)}_n)$, les coordonn\'ees homog\`enes d'un point  de l'espace
projectif $\P=\P (V^{r+1})=\P^{(r+1)(n+1)-1}$. Soit $\Sigma$ la sous-vari\'et\'e de $\P$ d\'efinie
comme le lieu des points 
$(x^{(0)},\ldots,x^{(r)})$ tels que les points $[x^{(0)}],\ldots,[x^{(r)}]$ de $\P V$ ne soient pas
lin\'eairement ind\'ependants. L'application
$$\rho:\P\moins\Sigma\lra G(r,\P V)$$
qui \`a $(x^{(0)},\ldots,x^{(r)})$ associe le $r$\tx plan engendr\'e par 
les points $[x^{(0)}],\ldots,[x^{(r)}]$ de $\P V$ est une fibration lisse connexe localement
triviale.

\ind Soient ${\bf f}=(f_1\ldots,f_s)$ les \'equations d\'efinissant $X$. On note $F$ l'adh\'erence
dans $\P$ de $\rho^{-1}(F_r(X))$; lorsque
$\delta(n,{\bf d},r)\ge 0$, il ressort du \theo\ \ref{fano} que la vari\'et\'e $F$ est
irr\'eductible de codimension ${{\bf d}+r\choose r}$ dans $\P$, lisse en dehors de $\Sigma$. 

\ind Pour tout entier $d$, on note ${\cal I}_d$ le sous ensemble de $\N^{r+1}$
form\'e  des
$(i_0,\ldots,i_r)$ tels que $\sum i_\a=d$; il a ${d+r\choose r}$ \'el\'ements. Pour tout \'el\'ement
$f$ de 
$\Sym^dV^*$ et tout \'el\'ement $I=(i_0,\ldots,i_r)$ de ${\cal I}_d$, on d\'efinit un
polyn\^ome $f_I$ multihomog\`ene de mutidegr\'e $(i_0,\ldots,i_r)$ sur $\P$ en posant
$$f(\l_0x^{(0)}+\cdots+\l_rx^{(r)})=\sum_{I\in {\cal I}_d}\l^If_I(x^{(0)},\ldots,x^{(r)})\ ,
\leqno{\formule}\hbox{\label{deffI}}$$
o\`u $\l^I=\l_0^{i_0}\cdots\l_r^{i_r}$; on convient aussi que $f_I$ est nul si l'un des $i_\a$ est
strictement n\'egatif. En dehors de $\Sigma$, la vari\'et\'e $F$  est d\'efinie par les \'equations
$$f_i(\l_0x^{(0)}+\cdots+\l_rx^{(r)})=0\qquad {\rm pour}\quad i=1,\ldots,s\quad{\rm et\ pour\
tout}\quad (\l_0,\ldots,\l_r)\in\P^r\ ,$$
\cad\ par les ${{\bf d}+r\choose r}$ \'equations $f_{i,I}$, pour $i=1,\ldots,s$ et $I\in 
{\cal I}_{d_i }$. En fait, comme $\Sigma$ est de codimension $n-r$ dans $\P$, si on suppose
$n-r>{{\bf d}+r\choose r}$, ces \'equations d\'efinissent $F$ dans $\P$; la vari\'et\'e $F$ est
alors une intersection compl\`ete irr\'eductible, lisse en dehors de $\Sigma$.
Son
degr\'e est la suite
${\bf D}$ o\`u chaque
$d_i$ est r\'ep\'et\'e ${d_i +r\choose r}$ fois. Posons $r_1=(r({\bf D})+1)(r+1)-1$; on suppose $ 
\d (n,{\bf d},r_1)
\ge 0$, de sorte qu'il existe un $r_1$\tx plan
$\Lambda_1=\P W_1$ contenu dans
$X$; on le suppose d\'efini par les \'equations
$x_{r_1+1}=\cdots =x_n=0$. On note $\Lambda^{r+1}_1$ le $((r_1+1)(r+1)-1)$\tx plan $\P (W_1^{r+1})$
dans
$\P$. 

 \ind Soit 
$\Lambda   $ un
$r({\bf D})$\tx plan contenu dans
$\Lambda^{r+1}_1$ et disjoint de $\Sigma$ (on pr\'ecisera plus bas notre choix de $\Lambda$). 
En vue d'appliquer le \theo\ \ref{PS}, on veut v\'erifier que l'application
$\beta:k^{(r+1)(n+1)}\ra \Gamma_\Lambda({\bf D}-1) $ d\'efinie
par
$$\beta(a^{(0)},\ldots,a^{(r)})
= \Bigl(
\sum_{j,l} a_l^{(j)}
\Bigl({\partial f_{i,I}\over\partial z_l^{(j)}}\Bigr)
_{\displaystyle{\vert_\Lambda}}\Bigr)_{1\le i\le s,\ I\in{\cal
I}_{d_i }}
$$ est surjective. D\'erivons l'\'egalit\'e \ref{deffI} par rapport \`a
$x^{(j)}_l$; on obtient
$$\l_j{\partial f\over\partial z_l}(\l_0x^{(0)}+\cdots+\l_rx^{(r)})=\sum_{I\in {\cal
I}_d}\l^I {\partial f_I\over\partial z^{(j)}_l}(x^{(0)},\ldots,x^{(r)})\ ,$$ 
de sorte que si  $\eps_j$ est l'\'el\'ement de ${\cal I}_1$ dont toutes les
composantes sauf la $j$i\`eme sont nulles, on a
$$\Bigl({\partial f\over\partial z_l}\Bigr)_{I-\eps_j}={\partial f_I\over\partial
z^{(j)}_l}\ ,$$ pour tout $I$ dans ${\cal
I}_d$ et tout $j=0,\ldots,r$. On peut donc \'ecrire
$$\beta(a^{(0)},\ldots,a^{(r)})
= 
\Bigl(
\sum_{j,l} a_l^{(j)}
 {\Bigl({\partial f_i\over\partial
z_l}\Bigr)_{I-\eps_j}}_{\displaystyle{\vert_\Lambda}}\Bigr)_{1\le i\le s,\  I\in{\cal I}_{d_i }}\ ,
$$ 
ou encore, en posant  $\partial_a f=
\sum_l a_l
\Bigl(\displaystyle{\partial f\over\partial z_l}\Bigr)
_{\displaystyle{\vert_{\Lambda_1}}}$ pour tout  
$f$ dans $\Sym^dV^*$, 
$$\beta(a^{(0)},\ldots,a^{(r)})
=  
\Bigl(
\sum_j {(\partial_{a^{(j)}}f_i)_{I-\eps_j}}_{\displaystyle{\vert_\Lambda}}\Bigr)_{1\le i\le
s,\  I\in{\cal I}_{d_i }} 
\ .$$

\th
Lemme
\enonce
Pour $n\ge r_1+{{\bf d}+r_1-1\choose r_1}$ et ${\bf f}$ g\'en\'erique dans $\Sym^{\bf d}V^*$ nul  sur
$\Lambda_1$, l'application
$$\eqalign{ \beta_1 :k^{ n+1  }  &\ \lra&  \Gamma_{\Lambda_1}({\bf d}-1)\ \  \cr
a\ \ \ &\ \longmapsto& ( \partial_af_1 ,\ldots, \partial_af_s )
\cr}$$ est surjective. 
\endth 

\ind  Il suffit de trouver un ${\bf f}$
pour lequel les  $\Bigl( \displaystyle{\partial f_1\over\partial z_l},\ldots,\displaystyle{\partial
f_s\over\partial z_l}\Bigr)_{0\le l\le n}$ engendrent $\Gamma_{\Lambda_1}({\bf d}-1 )$. 
Soient $J_1,\ldots, J_s$ des sous-ensembles disjoints de $\{ r_1+1,\ldots,n\}$ tels que $\Card
J_i= {d_i+r_1-1\choose r_1}$, et soit $\{ g_j\}_{j\in J_i}$ une base de $
\Gamma_{\Lambda_1}(d_i-1 )$. Il suffit de prendre $f_i=\sum_{j\in J_i} x_jg_j$.\cqfd

\ind Puisque $\Lambda$ est contenu dans $\Lambda_1^{r+1}$, l'application $\beta$ se factorise par
$(\beta_1)^{r+1}$, et il suffit de d\'emontrer que les applications
$$\eqalign{\gamma_d:\bigl( \Gamma_{\Lambda_1}(d-1)\bigr)^{r+1} &\ \lra & 
\Gamma_\Lambda(d-1)^{{\cal I}_d}\ \ \ \cr ( g^{(0)},\ldots,g^{(r)} )\ \ &\ \longmapsto&\Bigl(
\sum_j {g^{(j)}_{I-\eps_j}}_{\displaystyle{\vert_\Lambda}}\Bigr)_{I\in{\cal
I}_d}\cr}$$ sont surjectives pour $d=d_1,\ldots, d_s$. Nous allons montrer qu'elles sont surjectives
pour tout $d$, pour un choix convenable de $\Lambda$. Posons $x_{\a\b}=x_{\a (r({\bf D})+1)+\b}$, de
sorte que les $x_{\a\b}$, pour $0\le \a\le r$ et $0\le \b\le r({\bf D})$, forment des coordonn\'ees
sur $\Lambda_1$. Prenons pour $\Lambda$ le $r({\bf D})$\tx plan de $\Lambda_1^{r+1}$ d\'efini par les
\'equations  $$ x^{(j)}_{\a\b} =x^{(0)}_{0\b}\delta_{\a,j}\ ;$$
il est bien disjoint de $\Sigma$, et param\'etr\'e par les $y_\b=x^{(0)}_{0\b}$, pour
$\b=0,\ldots,r({\bf D}) $. 

\th
Lemme
\enonce
Pour tout entier  $d$, l'application
$$\eqalign{\gamma'_{d,q} :  \Gamma_{ \Lambda_1} (d-1)  &\ \lra &  \Gamma_\Lambda (d-1)^{{\cal
I}_{d-1}}\cr g &\ \longmapsto&\bigl(
    {g_I}_{\displaystyle{\vert_\Lambda}}\bigr)_{I\in{\cal
I}_{d-1}}\cr}$$ est surjective. 
\endth

\ind Soit $g=\displaystyle\prod_{\a,\b}
x_{\a\b}^{n_{\a\b}}$; on a 
$$g(\l_0x^{(0)}+\cdots+\l_rx^{(r)})_{\displaystyle{\vert_\Lambda}}=
\prod_{\a,\b}(\l_0x^{(0)}_{\a\b}+\cdots+\l_rx^{(r)}_{\a\b}
)^{n_{\a\b}}_{\displaystyle{\vert_\Lambda}}=
\prod_{\a,\b}(\l_\a y_\b 
)^{n_{\a\b}}\ ,$$
de sorte que  ${g_I}_{\displaystyle{\vert_\Lambda}}$ est le mon\^ome $\prod_\b y_\b 
 ^{\sum_\a n_{\a\b}}$ si $\sum_\b n_{\a\b}=i_\a$ pour tout $\a$, et est nul
sinon. Il reste \`a montrer que si $I=(i_0,\ldots,i_r)$ est fix\'e dans ${\cal I}_{d-1}$,  et si
$n_0,\ldots,n_{r({\bf D})}$ sont des entiers positifs de somme
$d-1$, il existe des entiers positifs $n_{\a\b}$ avec $\sum_\a n_{\a\b}= n_\b$ et $\sum_\b n_{\a\b}=
i_\a$ pour tous $\a$ et  $\b$, ce pour quoi 
il suffit de se donner deux partitions d'un ensemble \`a $n$ \'el\'ements 
en parties $(A_{\a})_{0\le\a\le r}$ et $(B_{\b})_{0\le\b\le r({\bf D})} $ de cardinaux respectifs
$i_{\a}$ et
$n_{\b}$,  et de prendre pour $n_{\a\b}$ le cardinal de $A_{\a}\cap B_{\b}$.\cqfd

\ind Pour montrer la surjectivit\'e de $\gamma_d$, il suffit donc de montrer celle de l'application
$$\eqalign{ \bigl( E^{{\cal
I}_{d-1}}\bigr)^{r+1}\ &\ \lra &  E^{{\cal
I}_d}\hskip 2cm\cr  
( {(g^{(0)}_I)}_I,\ldots,{(g^{(r)}_I)}_I ) &\ \longmapsto&\bigl(
    g^{(0)}_{J-\eps_0}+\cdots+ g^{(r)}_{J-\eps_r}
\bigr)_{J\in{\cal
I}_d}\cr}$$
o\`u $E$ est l'espace vectoriel $ \Gamma_{ \Lambda_1} (d-1)$; cela se fait sans difficult\'e pour
n'importe quel espace vectoriel $E$ par r\'ecurrence sur $r$.

 \ind On a montr\'e que toutes les
applications
$\gamma_{d_i}$, donc aussi l'application $\beta$, sont surjectives. Si $(r+1)(n+1)-1\ge n({\bf D})$,
on peut appliquer le \theo\ \ref{PS} pour conclure que
$F$ est unirationnelle, donc aussi $F_r(X)$; ceci termine  la d\'emonstration du
\theo .\cqfd

\vskip 1cm
\centerline{\bf R\'EF\'ERENCES}
\bigskip
\hangindent=1cm
[AK] A. Altman, S. Kleiman: {\sl Foundations of the theory
    of Fano schemes}, Comp. Math. {\bf 34} (1977), 3--47.
 
\hangindent=1cm
[BVV] W. Barth, A. Van de Ven: {\sl Fano varieties of lines 
    on hypersurfaces}, Arch. Math. {\bf 31} (1978), 96-104.


\hangindent=1cm
[BD] A. Beauville, R. Donagi: {\sl La vari\'et\'e des droites d'une hypersurface cubique de
dimension $4$}, C.R.A.S., t. 301, S\'erie I (1985), 703--706.

\hangindent=1cm
[BS] S. Bloch, V. Srinivas: {\sl Remarks on correspondences and algebraic cycles}, 
Am. J. of Math. {\bf 105} (1983), 1235--1253.

\hangindent=1cm
[BV] L. Bonavero, C. Voisin: {\sl Sch\'emas de Fano et
    vari\'et\'es de Moishezon}, C.R.A.S., \`a para\^\i tre.

\hangindent=1cm
[B1] C. Borcea: {\sl Deforming varieties of $k$\tx planes of projective complete intersections},
Pacific J. Math. {\bf 143} (1990), 25--36.

\hangindent=1cm
[B2] C. Borcea: {\sl Homogeneous Vector Bundles and Families of Calabi-Yau Threefolds. II}, in
Several Complex Variables and Complex Geometry (Santa Cruz 1989), Part 2, Proc.
Symp. Pure Math. {\bf 52} (1991), 83--91.

\hangindent=1cm
[Bo] R. Bott: {\sl Homogeneous vector bundles}, Ann. Math. {\bf 66} (1957), 203--248.

\hangindent=1cm
[C1] F. Campana: {\sl Remarques sur le rev\^etement universel des vari\'et\'es k\"ahleriennes
compactes}, Bull. S.M.F. {\bf 122} (1994), 255--284.

\hangindent=1cm
[C2] F. Campana: {\sl Connexit\'e rationnelle des vari\'et\'es de Fano}, Ann. Sci. E.N.S. {\bf 25}
(1992), 539--545.

\hangindent=1cm
[CTSSD] J.-L. Colliot-Th\'el\`ene, J.-J. Sansuc,  P. Swinnerton-Dyer: {\sl Intersection of two
quadrics and Ch\^atelet surfaces. I}, J.  reine  angew. Math. {\bf 373}
(1987), 37--107.

\hangindent=1cm
[D] O. Debarre: {\sl Th\'eor\`emes de connexit\'e pour les produits d'espaces projectifs et
les grassmanniennes}, Am. J. of Math. {\bf 118} (1996), \`a para\^\i tre.
  
\hangindent=1cm
[De] M. Demazure: {\sl A very simple proof of Bott's
    theorem}, Invent. Math. {\bf 33} (1976), 271--272.

\hangindent=1cm
[DR] U. Desale, S. Ramanan: {\sl Classification of Vector Bundles of Rank $2$ on Hyperelliptic
Curves}, Invent. Math. {\bf 38} (1976), 161--185.

\hangindent=1cm
[ELV] H. Esnault, M. Levine, E. Viehweg: {\sl Chow groups of projective varieties of very small
degrees}, Duke Math. J., \`a para\^itre.

\hangindent=1cm
[F] G.  Fano:  {\sl Sul sistema $\infty^2$ di rette contenuto in una variet\`a cubica generale
dello spazio a quattro dimensioni}, Atti Reale Accad. Sci. Torino {\bf 39} (1904), 778--792.

\hangindent=1cm
[Fu] W. Fulton: Intersection theory, Ergebnisse der Mathematik und
ihrer Grenzgebiete {\bf 2}, Springer Verlag, Berlin, 1984.

\hangindent=1cm
[FH] W. Fulton, J. Harris: Representation theory, Graduate Text in Mathematics {\bf
129}, Springer Verlag, Berlin, 1991. 

\hangindent=1cm
[H] D. Husemoller: Fibre bundles, 2nd ed., Graduate Text in Mathematics {\bf
20}, Springer Verlag, Berlin, 1966.

\hangindent=1cm
[K] J. Koll\'ar:  Rational Curves on Algebraic Varieties, Ergebnisse der Mathematik und
ihrer Grenzgebiete {\bf 32}, Springer Verlag, Berlin, 1996.

\hangindent=1cm
[KMM1] J. Koll\'ar, Y. Miyaoka, S. Mori: {\sl Rationally Connected Varieties}, J. Alg. Geom.
{\bf 1} (1992), 429--448.

\hangindent=1cm
[KMM2] J. Koll\'ar, Y. Miyaoka, S. Mori: {\sl Rational Connectedness and Boundedness of Fano
Manifolds}, J. Diff. Geom. {\bf 36} (1992), 765--769.

\hangindent=1cm
[L] Y. Laszlo: {\sl La dimension de l'espace des sections du diviseur th\^eta g\'en\'eralis\'e},
Bull. Soc. math. France {\bf 119} (1991), 293--306.

\hangindent=1cm
[M] I.G. Macdonald: Symmetric functions and Hall
  polynomials, Clarendon Press, Oxford, 1979.

\hangindent=1cm
[Ma1] L. Manivel: {\sl Th\'eor\`emes d'annulation pour les
    fibr\'es associ\'es \`a un fibr\'e ample}, Ann. Scuola Norm. Sup.  Pisa
  {\bf 19} (1992), 515--565.

\hangindent=1cm
[Ma2] L. Manivel: {\sl Applications de Gauss et pl\'ethysme
    2}, pr\'epublication de l'Institut Fourier {\bf 352}, 1996.

\hangindent=1cm [Mu] J. Murre: {\sl Discussion of a theorem of Morin},
notes de s\'eminaire \og Argomenti di Geometrica Algebrica\fg , Povo,
Trento, 1979.

\hangindent=1cm
[N] P.E. Newstead: {\sl Rationality of moduli spaces of stable bundles}, Math. Ann. {\bf 215}
(1975), 251--268.

\hangindent=1cm
[P] K. Paranjape: {\sl Cohomological and cycle-theoretic connectivity}, Ann. Math. {\bf 139}
(1994), 641--660.

\hangindent=1cm
[PS] K. Paranjape, V. Srinivas: {\sl Unirationality of the general Complete Intersection of
 small multidegree}, in Flips and Abundance
for Algebraic Threefolds, ed. J. Koll\'ar,
Ast\'erisque {\bf 211} (1992), 241--248.

\hangindent=1cm
[Pr] A. Predonzan: {\sl Intorno agli $S_k$ giacenti sulla variet\`a intersezione completa di pi\`u
forme}, Atti Accad. Naz. Lincei Rend. Cl. Sci. Fis. Mat. Natur.
  {\bf 5} (1948), 238--242. 

\hangindent=1cm
[R] S. Ramanan: {\sl The moduli spaces of vector bundles over an algebraic curve}, Math. Ann.
  {\bf 200} (1973), 69--84.

\hangindent=1cm
[S] A. Sommese: {\sl Complex Subspaces of Homogeneous Complex Manifolds II---Homotopy
Results}, Nagoya J. Math. {\bf 86}  (1982), 101--129.

\hangindent=1cm
[vW] B.L. van der Waerden: {\sl Zur algebraischen Geometrie
    2. Die geraden Linien auf den Hyperfl\"achen des $\P_n$}, Math.
  Ann. {\bf 108} (1933), 253--259.

\end